\documentclass[reqno]{amsart}
    \usepackage{float}
    \usepackage{amsmath}
    \usepackage{graphicx}
    \usepackage{latexsym}
    \usepackage{amsfonts}
    \usepackage{amssymb}
    \usepackage{verbatim} 
    \usepackage{color}  
    \theoremstyle{plain}
    \newtheorem{theorem}{Theorem}
    \newtheorem{corollary}[theorem]{Corollary}
    \newtheorem{lemma}[theorem]{Lemma}
    
    \theoremstyle{definition}

    \newtheorem{remark}[theorem]{Remark}
    \newtheorem*{remark*}{Remark}

    \newcommand{\pr}{\mathbf P}
    \newcommand{\e}{\mathbf E}
     
    \newcommand{\olF}{\overline F}

\begin{document}
    \title[Heavy traffic and heavy tails for the maximum]{Heavy
      traffic and heavy tails for 
      the maximum of a random walk}
    \author[Denisov]{Denis Denisov}
    \address{School of Mathematics, University of Manchester, Oxford Road, Manchester M13 9PL, UK}
    \email{denis.denisov@manchester.ac.uk}
    
    \author[Kugler]{Johannes Kugler}
    \address{Mathematical Institute, University of Munich, Theresienstrasse 39, D--80333
    Munich, Germany}
    \email{JKugler83@web.de}
    
    \begin{abstract}
    Consider a family of random walks $S_n^{(a)}=X_1^{(a)}+\cdots+X_n^{(a)}$ with negative drift
    $\mathbf E X_1^{(a)}=-a<0$ and finite variance 
    $\mbox{var}(X_1^{(a)})=\sigma^2<\infty$. 
    Let $M^{(a)}=\max_{n\ge 0} S_n^{(a)}$ be the maximums of the random walks. 
    The exponential asymptotics 
    $\pr(aM^{(a)}>x)\sim e^{-2x/\sigma^2}$, as $a\to 0$,  
    were found by Kingman and are known as 
    heavy traffic approximation in the queueing theory. 
    For  subexponential  random variables the large 
    deviation asymptotics for $\pr(M^{(a)}>x)\sim \frac{1}{a}\overline F^I(x)$  
    hold for fixed $a$ as $x\to\infty$.
    In this paper we present asymptotics for 
    $\pr(M^{(a)}>x)$,  which hold uniformly on the whole positive 
    axis, as $a\to 0$.  
    Thus, these uniform asymptotics include both the regime of normal and large deviations. 
    We identify the regions where exponential or subexponential asymptotics hold.   
    Our approach  is based on construction of corresponding super/sub - martingales 
    to obtain sharp upper and lower bounds. 
    \end{abstract}
    
    %\version\vspace{1cm}
    
    \maketitle
    {\scriptsize
    %TCIMACRO{\TeXButton{toc}{\tableofcontents}}%
    %BeginExpansion
    %\tableofcontents
    %EndExpansion
    }

    %
    %
    % TO DO:
    % 
    %
    %
    %
    %
    %
    %
    %
    %
    %
    %
    %
    %
    %
    %
    %
    %
    %
    %
    %
    
    %%%%%%%%%%%%%%%%%%%%%%%%%%%%%%%%%%%%%%%%%%%%%%%%%%%%%%%%%%%%%%%%%%%%%%%%%%%%%%%%%%%%%%%%%%%%%%%%%%%%%%%%%%
    
    \section{Introduction and statement of the results}
    \subsection{Introduction}
    Let $\{S_n^{(a)},n\geq0\}$ be a family of random walks with increments $X_i^{(a)}$ and starting point zero, that is,
    $$
     S_0^{(a)}:=0,\ S_n^{(a)}:=\sum_{i=1}^n X_i^{(a)},\ n\geq1.
    $$
    We shall assume that $X_1^{(a)},X_2^{(a)},\ldots$ are independent copies of a random variable 
    $X^{(a)} =X-a$ with a negative drift 
    $\mathbf{E}X^{(a)}=-a<0$, finite variance 
    $\sigma^2=\mbox{var} (X)$. We denote as $F$ the distribution 
    function of $X$. 
    For any $a>0$ the random walk $S_n^{(a)}$ 
    drifts to $-\infty$ almost surely  as $n\to \infty$, 
    and  
    therefore the  total maximum 
    $M^{(a)}:=\max_{k\geq0} S_k^{(a)}$ is finite almost surely. 
    However, $M^{(a)}\to \infty$ in probability as $a\to 0$.
    
    The distribution of the maximum of 
    a random walk appears in several classical models of applied probability, 
for example,  as the stationary distribution 
of the waiting time in a GI/GI/1 queue, or 
as the ruin probability in the Cram\'er-Lundberg risk process. 
For that reason it has attracted a lot of attention 
and the  exact and asymptotic behaviour of the distribution of the maximum 
has been extensively studied in the literature.

    The first asymptotic  result goes back, apparently, to Cram\'er and Lundberg (see,
    for example, \cite{A}).  
    Namely, if $a$ is fixed, 
    and there exists $h_0$ such that 
    \begin{equation}
    \label{Cramer}
    \mathbf{E}[e^{h_0X^{(a)}}]=1,
    \end{equation}
    and, in addition, $\mathbf{E}[X^{(a)}e^{h_0X^{(a)}}]<\infty$, then 
    \begin{equation}
    \label{Cramer-Lundberg}
    \mathbf{P}(M^{(a)}>x)\sim c_0e^{-h_0x}\quad\text{as }x\to\infty,
    \end{equation}
    where $c_0\in(0,1)$
    (here and throughout $a(x)\sim b(x), x\to \infty$ means that $a(x)/b(x)\to 1$ 
    as $x\to\infty$).

    If a solution to the equation \eqref{Cramer} 
    does not exist and, moreover, 
    $\mathbf{E}[\exp(\epsilon X)] = \infty$ for all
    $\epsilon>0$, 
    then, when $a$ is fixed,  
    the following asymptotics hold 
    \begin{equation}
    \label{vera}
      \mathbf{P}(M^{(a)}>x) \sim a^{-1}\olF^I(x) \quad \text{as } x \to \infty, 
    \end{equation}
    if the integrated tail of $F$, that is
    $\olF^I(x)= \int_x^\infty \olF(u) du$ with $\olF(x)=1-F(x)$, is subexponential, see 
    \cite{P75} or \cite{V77}. Here, we say that the distribution $F$ of a positive random variable is subexponential 
    if $\overline F^{*2}(x)\sim \overline F(x),$ as $x\to \infty$. 
    Main examples fo subexponential distributions 
    are Pareto ($\overline F(x)\sim x^{-r}, r>0$), lognormal and 
    Weibull ($\overline F(x)\sim e^{-x^\beta}, \beta\in(0,1)$). 
    We refer  to \cite{FKZ} for basic properties of the subexponential 
    distributions and for a proof of the subexponential asymptotics \eqref{vera}.  
    All subexponential distributions are heavy-tailed in the sense that  
    $\mathbf{E}[\exp(\epsilon X)] = \infty$ for all
    $\epsilon>0$. 
    
    Let us turn now to the case of $a\to 0$. This case is important and  interesting
    from the point of view of queueing theory as it describes the behaviour of a
    system in heavy traffic. There is a large volume of literature
    examining this case for a single-server queue and networks of queues. 
    These studies were initiated by 
    Kingman \cite{K61}, who considered the case when $|X^{(a)}|$ has an
    exponential moment and proved that 
    for fixed $x$, 
    \begin{equation}
    \label{asymp_heavy_traffic}
      \mathbf{P}(M^{(a)}>x/a) \sim e^{-2 x/\sigma^2} \quad \text{as } a \rightarrow 0.
    \end{equation}
    Kingman's approach was via the Wiener-Hopf factorisation, but he 
    also suggested an idea to derive this results from the corresponding 
    results for the maximum of the Brownian motion. The latter 
    approach was implemented by 
    Prohorov in \cite{P63} who applied the functional central limit theorem 
    to   extended \eqref{asymp_heavy_traffic} to the case that 
    the increments have finite variance, see \cite{SW11} for a recent discussion 
    of these two approaches.  
    Further ramifications to \eqref{asymp_heavy_traffic} 
    were obtained for an exponential family of distributions 
    by Borovkov in \cite[Theorem 3]{B64}, where 
    he obtained a complete asymptotic expansion 
    in $\theta$  for $\mathbf{P}(M^{(\theta)}>x/\theta)$ 
    using  the Wiener-Hopf factorisation. 
    In his work Borovkov used a different parametrisation 
    by the parameter $\theta$ of the exponential 
    family.  The first and the second term of the Borovkov expansion  
    were also obtained by considering renewal processes of ladder heights by Siegmund   
    in \cite[Theorem 2]{S79}. Recently, a complete asymptotic expansion for the exponential family 
    was obtained in \cite[Theorem 1]{BG06}.

    Asymptotics \eqref{asymp_heavy_traffic} can be obtained 
    by approximating the maximum of random walks with 
    the maximum of the Brownian motion with  a drift. 
    Thus, it can be viewed as a diffusion approximation for the maximum of the random walk. 
    On the other hand, asymptotics \eqref{vera} is explained 
    by a single large jump and  clearly represents the large deviations 
    regime for the maximum of the random walk. Thus, it is an interesting 
    question to identify when the normal deviations hold and 
    when the large deviations hold more precisely. 
    
    For the sums $S_n$ with zero mean and unit variance 
    the transition from normal approximation to large deviations has 
    been extensively studied. In particular, the  classical result of Nagaev 
    holds for the regularly varying distribution function $\overline F(x)\sim x^{-r}L(x), r>2 $ 
    and can be given as follows, 
    \begin{equation}
      \label{nagaev}
      \pr(S_n>x) \sim \overline \Phi(x/\sqrt n) + n\overline F(x), 
    \end{equation}  
    as $n\to \infty$ and $x/\sqrt n\to +\infty$, 
    see \cite{N63} or Theorem 1.9 in \cite{N79}.  
    Here  $\overline \Phi(t)=\frac{1}{\sqrt{2\pi}} \int_t^{+\infty} e^{-u^2/2}$ 
    it the tail of the distribution function of the standard normal law.  Thus, the first term 
    is \eqref{nagaev} corresponds to the Central Limit Theorem, while the second term represents 
    large deviations and is explained by a large jump in one of $n$ steps. 
    It is not difficult to identify the transition 
    point in \eqref{nagaev} and the regions of domination 
    of the normal approximations or large deviations approximations. 
    For lighter subexponential distributions, such as Weibull, the situation is more involved.  
    The general case of subexponential distributions is considered in 
    ~\cite{R93}, where an analogue of \eqref{nagaev} was obtained,  
     see also  \cite{DDS08} for a local version and further references. 

    For the maximum of the random walks  less is known about the 
    the transition from the normal approximation to large deviations, 
    that is from  the exponential asymptotics to subexponential asymptotics. 
    The {\em main objective}  
    of  this paper is to study transition phenomena 
    for the maximum of general random walks 
    for  a wide subclass of subexponential distributions, 
    including all standard examples.   
    In particular our  {\it main result} gives
     uniform in $x$ asymptotics for $\pr(M^{(a)}>x)$ 
     and can be viewed as an analogue of \eqref{nagaev} for the maximum 
     of the random walk. 
    
    The transition phenomena  are fairly well understood for the M/G/1 
    queues, see \cite{OBG}, \cite{OG}, \cite{BL} and references therein. 
    In \cite{OBG} and  \cite{OG} the case of the regularly varying distributions
    was considered.  In the case of regularly varying distribution  uniform asymptotics was obtained in 
    \cite[Theorem 2.2]{OBG} and the transition point from heavy traffic 
    regime to the large deviations was obtained in \cite[Corollary 2.3]{OBG}. 
        In \cite{BL} uniform asymptotics was obtained for 
    a general subclass of subexponential distributions including in particular 
    regular varying, lognormal and Weibull distributions. 
    The latter  uniform  asymptotics for the probability $\pr(M^{(a)}>x)$, 
    is a sum of the exponential term from  Kingman's asymptotics, 
    the integrated tail term  and a convolution term. 
    A more explicit representation together with the transition point 
     was obtained for regularly varying distributions  in Example 1.1 of \cite{BL}. 
     To derive the main result  of \cite{BL} Blanchet and Lam 
     applied the large deviation results for sums on the whole axis 
     from \cite{R93} and the results on the uniform renewal theory \cite{BG07}.

    For  the maximum of  general random walks 
    and, correspondingly, steady-state waiting time of $GI/GI/1 $ queues, 
    the transition  phenomena have been studied only for regularly varying distributions. 
    In particular, Theorem~4 of  
     \cite{KW} states that     
     if the increments are regularly varying of index $r>2$, that is,
    $$
     \pr(X>u)= u^{-r}L(u)
    $$
    where $L$ is a  slowly varying function, then 
    $$
    \pr(M^{(a)}>x_a) \sim \frac{1}{a} \overline F^{I}(x_a)
    $$
    for $x_a$ such that 
    $$
    \liminf_{x\to \infty} \frac{x_a}{a^{-1}\ln{1/a}}>e^r\frac{r-2}{2}\sigma^2 
    $$
    extending thus asymptotics \eqref{vera} to a larger region. 
    The methodology of \cite{KW} was based on obtaining analogues 
    of Fuk-Nagaev's inequality for maximums and then passing to the limit. 
    In Theorems 3.9 and 3.11 of 
    \cite{K14} methods similar to \cite{KW} were used to establish 
    that 
    $$
    \pr(M^{(a)}>x_a) \sim \frac{1}{a} \overline F^{I}(x_a)
    $$
    for $x_a$ such that 
    $$
    \liminf_{x\to \infty} \frac{x_a}{a^{-1}\ln{1/a}}>\frac{r-2}{2}\sigma^2 
    $$
    and 
    $$
    \pr(M^{(a)}>x_a) \sim e^{-2ax_a/\sigma^2}
    $$
    for $x_a$ such that 
    $$
    \limsup_{x\to \infty} \frac{x_a}{a^{-1}\ln{1/a}}<\frac{r-2}{2}\sigma^2 
    $$
    establishing, in particular, that for regularly varying distributions 
    $$
    x_{RV}(a) = \frac{r-2}{2}\sigma^2  a^{-1}\ln(1/a)
    $$
    is a transition point between the heavy  traffic regime and the subexponential regime.

    Before stating our main results we will briefly discuss why the case of  M/G/1 queue 
    has been mainly considered.  
    This is due to the existence of  the representation of $M^{(a)}$ 
    as a geometric sum of independent and identically distributed random variables:
    \begin{equation}
      \label{eq.geom.sum}
     \pr(M^{(p)}>x) = \sum_{k=0}^\infty p(1-p)^k \pr(\chi_1^+ + \chi_2^+ \dots + \chi_k^+ >x),
    \end{equation}
    where $\chi_1$ is the ascending ladder height and $p=\pr(M=0)$.  
    In the case of M/G/1 queue the distribution 
    of $\chi_1$ is known explicitly and is proportional to the integrated tail $\overline F^I$. 
    Thus it is sufficient to analyse 
    the  distribution of the geometric sum, when $p\to 0.$ 
    In the general case this approach is still possible but seems to be rather complicated 
    due to the fact that the distribution of $\chi_1^+$ is  less explicit.  
    Thus we use use a different approach which relies on
    martingale methods.

    \subsection{Statement of main results}
    We will first present a class of subexponential distributions that we will consider. 
    Let 
    \begin{equation}
      \label{sc1}
    \overline F(x) \sim e^{-g(x)}x^{-2}, \quad  x\to \infty,
    \end{equation}
    where $g(x)$ is a continuously differentiable eventually increasing function such that 
    eventually 
    \begin{equation}
      \label{sc2}
    \frac{g(x)}{x^{\gamma_0}}\downarrow 0, \quad x\to\infty,
    \end{equation}
    for some $\gamma_0\in(0,1)$. 
  Clearly, monotonicity in \eqref{sc2} implies 
    \begin{equation}
      \label{sc3}
      g'(x)\le \gamma_0 \frac{g(x)}{x}
    \end{equation}  
    for all sufficiently large $x$. 
    Due to the asymptotic nature of equivalence in \eqref{sc1} 
    one  may assume without loss of generality that 
    \eqref{sc2} and \eqref{sc3} hold for all $x>0$. 
    Using the Karamata representation theorem one can show that 
    this class of subexponential distributions  include regularly varying  
    distributions $\overline F(x)\sim x^{-r}L(x),$ for $r>2$. Also, it is not difficult 
    to show that lognormal distributions and  Weibull distributions 
    ($\overline F(x) \sim e^{-x^\beta},\beta\in(0,1)$) 
    belong to our class of distributions. 
    Previously this class appeared in \cite{R93} 
    for the analysis of 
    large deviations of sums of subexponential random variables on the whole axis. 

    Further, for each $a>0$,  
    let $x(a)$ be a solution to 
    \begin{equation}
      \label{eq.boundary}
    \frac{2a}{\sigma^2} = \frac{g(x)}{x}
    \end{equation}
    Since, by \eqref{sc2} the function $g(x)/x$ is monotone decreasing to $0$ 
    we can pick $x(a)$ in such a way that $x(a)\uparrow \infty$ 
    as $a\downarrow 0.$  (here and throughout $f(a)\uparrow a_0$ 
    means that $f$ is monotonically increasing to $a_0$ and 
    correspondingly  $f(a)\downarrow a_0$ 
    means that $f$ is monotone decreasing to $a_0$). 
    We will show that $x(a)$ is the transition point 
    between the heavy traffic regime and the subexponential regime. 

    We will also need the following version of 
    the Cram\'er equation \eqref{Cramer} for a truncated random variable.
    Let $\theta_a$ be a function such that 
    \begin{equation}
      \label{def.theta.xa}
       \mathbf E\left[e^{\theta_a X^{(a)}} ;X^{(a)}\le 1/a\right]=1+o(a/x(a))
       ,\quad a\to 0. 
      \end{equation} 
    It is not difficult to show that even the exact solution 
    $$
    \mathbf E\left[e^{\theta^{exact}_a X^{(a)}} ;X^{(a)}\le 1/a\right]=1
    $$  
    exists. However, the form \eqref{def.theta.xa} is more convenient, 
    as this approximate equation is easier to   
    solve than the exact equation. 
    Using the Taylor expansion it is not difficult to show that 
    $\theta_a\sim 2a/\sigma^2,$ as $a\to 0$. 
    In Section~\ref{sec.sol.theta.xa} we will obtain a polynomial (in $a$) 
    solution to the  equation \eqref{def.theta.xa} with sufficient for our purposes precision. 

    %Assuming  existence of sufficiently many moments 
    %we can show that 
    %one can pick a
    %%solution $\theta_a$ of the form
    %$$   
    %\theta_a=\frac{2a}{\sigma^2}+C_2a^2+\cdots+C_{k^*}a^{k^*} +o(a^{k*}). 
    %$$
    %In more details the equation \eqref{def.theta.xa} will be considered 
    %in Section~\ref{sec.sol.theta.xa}.    

    We are now in position to state our first result. 
    \begin{theorem}
      \label{thm.main}
      Let $\e[X]=0$, $\mbox{var}(X)=\sigma^2$ and 
      $\e[(X^+)^{2+\varepsilon}]$ for some $\varepsilon>0$. 
      Let the distribution function 
      $F$ satisfy condition \eqref{sc1} and \eqref{sc2}.  
      Let $x(a)$  and $\theta_a$ be the functions defined above 
      by \eqref{eq.boundary} and \eqref{def.theta.xa} respectively.  
      Then, for any $\delta\in (0,1)$ 
      and $A_a\uparrow\infty$, uniformly 
      in $x\ge 0$, as $a\to 0$,  
      \begin{align}
        \nonumber
        \pr(M^{(a)}>x) 
        &\sim e^{-\theta_a x}\mathbf{1}\left\{x<(1+\delta)x(a)\right\} \\
        \label{eq.main}
        &+\frac{1}{a} \overline F^I(x)\mathbf{1}\left\{x>(1-\delta)x(a)\right\}\\ 
        &+
        \frac{2}{\sigma^2}\int_{(1-\delta)x(a)}^x 
        \left(
          \frac{1}{\theta_a}+(x-y)
        \right)
        \overline F(y)dy 
        \mathbf{1}\left\{(1-\delta)x(a)\le x\le A_ax(a)\right\}. 
        \nonumber
      \end{align}  
    \end{theorem}  

    \begin{remark}
      \label{R1}
       In \cite{BL} a similar result was proved for the geometric sum 
       \eqref{eq.geom.sum} (and thus of M/G/1 queue). 
       This result requires to solve an analogue of 
       of \eqref{def.theta.xa} for the increasing ladder height 
       $\e[e^{\theta^*_p \chi^+_1}; \chi^+_1\le 1/p]=1/(1-p)$. 
       Then, assuming $2+\varepsilon$ moments of $\chi^+$ and under similar 
       assumption on $\overline F$, they show that 
       uniformly in $x$,  as $p\to 0$,
       \begin{align*}
       \pr(M^{(p)}>x) 
       &\sim e^{-\theta_p^* x} \\
       &+ 
        \left(
        \frac{\olF^I(x)}{p}+  
        \int_{1/p}^x
        \left( 
        \frac{1}{p}+\frac{x-y}{\e\chi_1^+}\right) e^{-\theta_p^*(x-y)} 
        d\pr(\chi^+\in dy)
        \right)\mathbf{1}\left\{x\ge \frac{1}{p}\right\}.
      \end{align*}
      Note that the parametrisation is slightly different ($p\to 0$). 
      The latter asymptotics was then analysed for the case of regularly varying $F$. 
       \end{remark}

       Asymptotics \eqref{eq.main} contains three terms: 
       the exponential term, corresponding to the Kingman asymptotics, 
       subexponential term  and a  convolution term. 

       Under further assumptions we can give a more explicit form  
       of \eqref{eq.main}. 
       First, we will study the case when the convolution  
       term disappears. 
       \begin{corollary}
        \label{thm.rv}
        Let $\e[X]=0$, $\mbox{var}(X)=\sigma^2$ and 
        $\e[|X|^{2+\varepsilon}]$ for some $\varepsilon>0$.  
        Let $F$ satisfy \eqref{sc1}, where 
        $g(x)$ is continuously differentiable and 
        $$
        g'(x) = o(g(x)/x),\quad x\to \infty.
        $$
        Let $x(a)$ be defined by \eqref{eq.boundary}. 
        Then, for any $\delta \in (0,1)$, 
        uniformly in x, as $a\to 0$, 
        \begin{align}
          \label{eq.main.rv}
          \pr(M^{(a)}>x) 
          &\sim e^{- 2ax/\sigma^2}\mathbf{1}(x<(1+\delta)x(a)) 
            +\frac{1}{a} \overline F^I(x)\mathbf{1}(x>(1-\delta)x(a)). 
        \end{align} 
       \end{corollary} 
       This theorem includes regularly varying distributions 
       and lognormal distributions. 
       For Weibull distributions 
       the convolution term plays a role, as explained 
       by the following statement  
       \begin{corollary}
        \label{thm.weibull}
        Let the conditions of Theorem~\ref{thm.main} hold. 
        Assume in addition that  
        $$
        g'(x) = x^{\beta-1}L(x),\quad x\to \infty,
        $$
        where $L(x)$ is a slowly varying function 
        and $\beta\in (0,1)$. 
        Then, for any $\delta \in (0,1-\beta)$ and 
        any function $A_a\uparrow\infty$, 
        uniformly in x, as $a\to 0$, 
        \begin{align}
          \label{eq.main.weibull}
          &\pr(M^{(a)}>x) 
          \sim e^{-\theta_a x}\mathbf{1}\left\{x<(1+\delta)x(a)\right\} 
            +\frac{1}{a} \overline F^I(x)\mathbf{1}\left\{x>(1-\delta)x(a)\right\}\\ 
           &+ 
           \frac{2}{\sigma^2}
           \overline F^{I}(x) g'(x)
 \left(
  \frac{1}{\theta_a(\theta_a-\beta\frac{g(x)}{x})}
  +
 \frac{1}{(\theta_a-\beta\frac{g(x)}{x})^2}
 \right)
 \mathbf{1}\left\{(1-\delta)x(a)<x<A_ax(a)\right\}.
          \nonumber
        \end{align} 
       \end{corollary} 
    \begin{remark}
      It is clear from \eqref{eq.boundary} that 
      in Corollary~\ref{thm.weibull} the boundary 
      $$
      x(a) \sim 
      \left(
        \frac{1}{a}
        \right)^{1/(1-\beta)} L^*(1/a), 
        \quad a\to 0, 
      $$  
      where $L^*(\cdot)$ is a slowly  varying function. 
      Then, if $\beta<1/2$ and 
      $\e[|X|^{\frac{2-\beta}{1-\beta}+\varepsilon}]<\infty$
      for some $\varepsilon>0$ then it follows 
      from Lemma~\ref{lem.eq.weibull}  
      that one can take in Corollary~\ref{thm.weibull} 
      $
      \theta_a= 2a/\sigma^2+\sum_{j=2}^{[\frac{1}{1-\beta}]}C_ja^j
      $, 
      where $[x]$ is the integer part of $x$.   
      In particular, for $\beta<1/2$ one can put $\theta_a=2a/\sigma^2$. 
      \end{remark}   

     \subsection{Sub/super-martingale construction}
    In this subsection we formulate our most general theorem, from which we will 
    derive Theorem~\ref{thm.main}. 
    The proof of this general theorem 
    will be  done via a construction 
    of sub/super martingales which approximate our asymptotics.    
    
    First we need some notation. 
    Recall that  $F$ is  the distribution function of $X$ and let 
    $F_a$ be the distribution function of $X^{(a)}$. 
    Clearly, for any $x$, we have the equality $F_a(x) = F(x+a).$  
    Let $\overline F(x) = 1-F(x)$ and $\overline F_a(x) = 1-F_a(x)$.
    Let 
    $\varepsilon_a$ 
    be a function such that $\varepsilon_a\downarrow 0,$ 
    as $a\to 0$ 
    and let 
    $$
    \overline F_+(x) = \overline F\left(
      x+\frac{\varepsilon_a}{a}
      \right),\quad 
      \overline F_+^I(x) = \int_x^\infty \overline F_+(y) dy. 
      \quad  x\ge 0. 
    $$ 
    Let $c_a$ be a function such that $c_a=o(a)$ and 
    let $\theta_a$ be a function such that 
    \begin{equation}
    \label{def.theta}
     \mathbf E\left[e^{\theta_a X^{(a)}} ;X^{(a)}\le 1/a\right]=1+o(ac_a)
     ,\quad a\to 0. 
    \end{equation} 
     One can easily see by Taylor expansion, 
     that $\theta_a\sim 2a/\sigma^2$ as $a\to 0$. 

     We will assume that the distribution function $F$ 
     satisfies 
     \begin{equation}
      \label{subex}
      \sup_{x\ge 2\varepsilon_a/a}
      \frac{\int_{\varepsilon_a/a}^{x-\varepsilon_a/a} dy \overline F(y)\overline F(x-y)}
      {\overline F(x)} = o(c_a).
    \end{equation} 
    \begin{remark}
    Equation \eqref{subex} defines a subclass of $\mathcal S^*$ 
    in a way similar to the introduction 
    of a subclass of $\mathcal S$ in equation (4) of \cite{DDS08}. 
    $\mathcal S^*$ is a subclass of subexponential 
    distributions   $\mathcal S$  introduced by 
    Kl\"uppelberg~\cite{Kl88}. It is known that it is a proper subclass of 
    $\mathcal S$, that is there exist distributions with finite mean 
    that belong to $\mathcal S$, but do not belong to $\mathcal S^*$, see 
    Section~6 of \cite{DFK}. 
    Class of distributions satisfying \eqref{subex} is sufficiently 
    rich and includes all   major subexponential distributions, see 
    Lemma~\ref{lem.sequences} below. 
    \end{remark}  

     Next let $\alpha>0$ be a fixed number. 
   Let 
     $\lambda^\pm_a = \theta_a\pm c_a$
    and define  functions $\overline G_+$ and $\overline G_-$ as follows  
    \begin{equation}
      \label{def_G}
      \overline G_\pm(x)
      =\begin{cases}
          e^{\pm 2\alpha},& x<0, \\ 
          e^{\pm\alpha}e^{-\lambda^\pm_a x }
          +\frac{2}{\sigma^2\lambda^\pm_a} 
          \overline F_+^I(x)
          +\frac{2}{\sigma^2}\int_0^x dz e^{-\lambda^\pm_a (x-z)}
          \left(
            \frac{1}{\lambda^\pm_a } +(x-z) 
          \right)
          \overline F_+(z), & x\ge 0, 
        \end{cases} 
    \end{equation} 
    where for $\overline G_+$ we use 
    $e^{+2\alpha}, e^{+\alpha}$ and 
    $\lambda_+$ and 
    for $\overline G_-$ we use 
    $e^{-2\alpha}, e^{-\alpha}$ and 
    $\lambda_-$. We use the same functions 
    $\overline F_+$ and $\overline F^I_+$ for both $\overline G_+$ 
    and $\overline G_-$.

  \begin{theorem}
    \label{prop1}
    Let $\e[X]=0$ and $\mbox{var} (X)=\sigma^2$. 
    Let  $c_a$ and $\varepsilon_a$ be functions such that 
    $c_a=o(a)$ and  $\varepsilon_a\downarrow  0$.  
    Let $\theta_a$ be a function 
    satisfying \eqref{def.theta} 
    such that $\theta_a\sim 2a/\sigma^2,$ as $a\to 0$. 
    Assume that  
    $\overline F(\varepsilon_a/a) = o(ac_a)$ 
    and $\overline F^I(\varepsilon_a/a)=o(a)$, 
    as $a\to 0$.  
    Then, for any 
    $\alpha>0$,  there exist $a_0>0$  such that 
    for any $a\in(0,a_0)$, 
    \begin{equation}
      \label{eq.lower.bound}
      \pr(M^{(a)}>x) \ge e^{-2\alpha}\overline G_+(x+\varepsilon_a/a),  
    \end{equation}   
    and, if we assume in addition the \eqref{subex} holds 
     then 
    \begin{equation}
      \label{eq.upper.bound}
      \pr(M^{(a)}>x) \le e^{2\alpha}\overline G_-(x+\varepsilon_a/a),  
    \end{equation}   
    where $\overline G_\pm$ are defined according to 
    \eqref{def_G} with the above $c_a,\theta_a,\varepsilon_a$ and $\alpha$. 
  \end{theorem}  
  The accuracy of this theorem depends, of course, 
  on a good choice of functions $c_a$ and $\varepsilon_a$. 
In the proof of Theorem~\ref{thm.main} we demonstrate how these function can be picked in order 
to obtain a more explicit statement \eqref{eq.main}. 

\begin{comment}

  In what follows we will demonstrate that 
  for sufficiently regular distributions 
  we can define $c_a$ and $\varepsilon_a$ 
  in such a way that 
  the upper and lower bounds become asymptotically equivalent 
  and 
  \begin{equation}
    \label{asymp}
    \pr(M^{(a)}>x) \sim 
    e^{-\theta_a x}
    +
    \frac{2}{\sigma^2 \theta_a} 
          \overline F^I(x)
          +\frac{2}{\sigma^2}\int_0^x dz e^{-\theta_a (x-z)}
          \left(
            \frac{1}{\theta_a } +(x-z) 
          \right)
          \overline F(z)dz
  \end{equation}  
  uniformly in $x>0$ as $a\to0$.
\end{comment}
  
    \subsection{Further discussion}
    The approach presented in this paper 
    is based on construction of supermartingales and submartingales 
    in order to obtain accurate lower and upper bounds. 
    Previously, this approach was used by one of the authors 
    and V. Wachtel in \cite{DW12} to prove \eqref{vera}, 
    which gives subexponential asymptotics for $\pr(M^{(a)}>x)$   
    when $a$ is fixed. 
    For diffusion approximation this approach works as well. 
    Indeed, if one considers Brownian motion with a negative drift 
    $X_t=\sigma B_t-at$, then by considering an exponential martingale 
    and applying  Doob's optional theorem it is straightforward to show 
    that 
    $$
    \pr(\sup X_t >x) = e^{-2ax/\sigma^2}, \quad x\ge 0.
    $$
    One can immediately recongnise the Kingman asymptotics on the right-hand-side. 
    In the case of random walks we cannot apply an exponential martingale immediately 
    since we have to deal with overshoots, but we can approximate this martingale 
    with appropriate super- and submartingales.   

    In general, for diffusion processes it is natural to use martingales in order to obtain 
    information about exit times and probabilities of return from  various domains. 
    Thus, we can use   a similar technique to estimate return probabilities  
    for Markov chains. However, in order to deal with overshoots we should 
    use super/sub martingales constructed on the base of corresponding diffusion martingale. 
    This idea can be used to approximate 
    return probabilities of Lamperti Markov chains and to prove renewal theorems, 
    as in \cite{DKW13} and \cite{DKW16}, show transience of Markov chains, 
    see \cite{DKW16} and to some extent \cite{DF},  \cite{FD}.  

    There is a traditional approach to the analysis of the maximum 
    of a random walk via an analysis of geometric random sums. 
    This approach seems to be very appealing to give a proof of Theorem~\ref{thm.main}, 
    as a similar theorem was obtained  in \cite{BL} for geometric sums. 
    To apply this approach one  first needs to find a relationship between 
    the extensions of Cram\'er equation \eqref{def.theta.xa} and the 
    corresponding equation for the ascending ladder height $\chi^+_1$. 
    It seems to be likely that this relationship exists and one can show that 
    if $\theta_a$ solves \eqref{def.theta.xa} then the same $\theta_a$ 
    will solve 
    $$
    \e[e^{\theta_a\chi_1^+};\chi_1^+\le 1/a] = 1+ o(a/x(a)).
    $$
    The latter equation together with a uniform  (in $a$) 
    control of the renewal 
    function  of descending ladder heights 
    might give an alternative derivation of Theorem~\ref{thm.main}. 

    In the proof of the most general theorem Theorem~\ref{prop1} 
    we used a uniform version of the class $\mathcal S^*$. 
    It is known that if $F\in \mathcal S^*$ then 
    for non-lattice distributions  
    $$
    \pr(M^{(a)}\in (x,x+T]) \sim \frac{T}{a}\overline F(x), \quad 
    x\to \infty, 
    $$
    for  any fixed $T$. Thus, when $a$ is fixed, the class $\mathcal S^*$ gives a local, 
    rather  than global asymptotics for distribution of $M$. 
    Hence, it seems reasonable to expect that the main results of the present 
    could be strengthened to give uniform (in $x$) asymptotics  
    for $\pr(M^{(a)}\in (x,x+T])$, as $a\to 0$. Local asymptotics requires 
    more precise arguments and will be considered elsewhere.

    %%%%%%%%%%%%%%%%%%%%%%%%%%%%%%%%%%%%%%%%%%%%%%%%%%%%%%%%%%%%%%%%%%%%%%%%%%%%%%%%%%%%%%%%%%%%%%%%%%%%%%%%%%%%
    
    \subsection{Organisation of the paper}
    The paper is  organised as follows. 
    First we present the sub/super-martingale construction and 
    prove the most general Theorem~\ref{prop1} in Section~\ref{sec.prop1}.  
    Next, in Section~\ref{sec.main} 
     we prove Theorem~\ref{thm.main} by 
     specialising  Theorem~\ref{prop1} to the subclass of subexponential distributions satisfying 
    \eqref{sc1} and \eqref{sc2}.  
    After that, in Section~\ref{sec.sol.theta.xa} 
    we study solutions to the equation~\eqref{def.theta.xa}. 
    Finally, in Section~\ref{sec.final} 
    we prove Corollary~\ref{thm.rv} and Corollary~\ref{thm.weibull}.

    \section{Proof of Theorem~\ref{prop1}}
    \label{sec.prop1}

    \subsection{Perturbation of $\theta_a$}
    In order to obtain super/sub martingales 
    we are planning to slightly increase/decrease the parameter  
    $\theta_a$ 
    and the next lemma  shows  that the equation \eqref{def.theta} 
    will still be approximately correct. 
    \begin{lemma}
      \label{lem7}
      Let $\theta_a$ be defined according 
      to \eqref{def.theta} and $c_a=o(a)$. 
        Then, 
        \begin{align}
          \label{eq.theta.below2}
            \e\left[e^{(\theta_a+c_a) X^{(a)}};X^{(a)}\le 1/a\right] &\ge     1 + a c_a + o(a c_a), \quad a\to 0,\\  
            \e\left[e^{(\theta_a-c_a) X^{(a)}};X^{(a)}\le 1/a\right] &\le     1 - a c_a + o(a c_a), \quad a\to 0
            \label{eq.theta.above2}
        \end{align} 
    \end{lemma}  
    \begin{proof}
      We will first derive the lower bound \eqref{eq.theta.below2}. 
      Using the elementary inequality $e^x\ge 1+x$ 
      two times and \eqref{def.theta} we obtain,     
    \begin{align*}
        &\e\left[e^{(\theta_a+c_a) X^{(a)}};X^{(a)}\le 1/a\right]
        \ge \e\left[e^{\theta_a X^{(a)}};X^{(a)}\le 1/a\right] 
         + c_a \e\left[X^{(a)}e^{\theta_a X^{(a)}};X^{(a)}\le 1/a\right]\\
        &\hspace{1cm}\ge 1 + o(ac_a)+
        c_a \e\left[X^{(a)};X^{(a)}\le 1/a\right] + \theta_a c_a \e\left[(X^{(a)})^2;X^{(a)}\le 1/a\right].
      \end{align*} 
    Now note that  since 
    the family of random variables $\{(X^{(a)})^2\}_{a>0}$ 
    is uniformly integrable, 
    \begin{align}
     \nonumber 
    \e\left[X^{(a)};X^{(a)}\le 1/a\right] 
    &=\e\left[X^{(a)}\right]-\e\left[X^{(a)};X^{(a)}> 1/a\right]\\
    &\ge -a- a\e\left[(X^{(a)})^2;X^{(a)}> 1/a\right]
    =-a-o(a)
    \label{eq.first}
    \end{align}
    and 
    \begin{equation}
    \e\left[(X^{(a)})^2;X^{(a)}\le 1/a\right] 
    =\sigma^2 +a^2  
    -\e\left[(X^{(a)})^2;X^{(a)}> 1/a\right] 
    =\sigma^2+o(1).
    \label{eq.second}
    \end{equation}
    Hence, recalling that $\theta_a \sim \frac{2a}{\sigma^2}$ we obtain 
    $$
    \e\left[e^{(\theta_a+c_a) X^{(a)}};X^{(a)}\le 1/a\right] \ge 
    1 - a c_a + \theta_a c_a \sigma^2 + o(a c_a)
    =1 + a c_a + o(a c_a),
    $$
    proving thus \eqref{eq.theta.below2}. 

    To prove the upper bound \eqref{eq.theta.above2} 
    we will use the elementary estimate $e^x\le 1+x+x^2$, which is valid for $x\le 1$
    and then again the estimate $e^x\ge 1+x$, which is always valid. We have, 
    \begin{align*}
      &\e\left[e^{(\theta_a-c_a) X^{(a)}};X^{(a)}\le 1/a\right] 
      \le \e\left[e^{\theta_a X^{(a)}};X^{(a)}\le 1/a\right] 
       - c_a \e\left[X^{(a)}e^{\theta_a X^{(a)}};X^{(a)}\le 1/a\right]\\
      &\hspace{1cm}+ c_a^2 \e\left[(X^{(a)})^2e^{\theta_a X^{(a)}};X^{(a)}\le 1/a\right] \\  
      &\hspace{0.5cm}\le 1 +o(ac_a)
      - c_a \e\left[X^{(a)};X^{(a)}\le 1/a\right] 
      - \theta_a c_a \e\left[(X^{(a)})^2;X^{(a)}\le 1/a\right]\\
      &\hspace{1cm}+ c_a^2 e^{\theta_a/a}\e\left[(X^{(a)})^2\right].
    \end{align*} 
    Recalling that $\theta_a \sim \frac{2a}{\sigma^2}$ and using 
    \eqref{eq.first} and \eqref{eq.second} we obtain
    $$
    \e\left[e^{(\theta_a-c_a) X^{(a)}};X^{(a)}\le 1/a\right] 
    \le 
    1 + a c_a - \theta_a c_a \sigma^2 + o(a c_a)
    =1 - a c_a + o(a c_a),
    $$
    as required. 
  \end{proof}  
  Let 
    \begin{equation}
      \label{def_f_star}
      F_a^* (t)   = \begin{cases}
                      \overline F_a(t),& t\ge 0 \\ 
                      -F_a(t),& t<0. 
                    \end{cases}
    \end{equation}  
  Function  $F_a^* (\cdot)$ 
  appears after the integration by parts of $\overline G_\pm$. 
  Next Lemma  reformulates the  equation \eqref{def.theta} 
  in terms of $F_a^* (\cdot)$ and follows from 
  Lemma~\ref{lem7}. 
  \begin{lemma}
      Let  $\lambda_a^+=\theta_a+c_a$ and 
      $\lambda_a^-=\theta_a-c_a$. 
        Assume that $c_a=o(a)$, $\varepsilon_a\downarrow 0$ 
        and $\overline F(\varepsilon_a/a)=o(ac_a)$.  
        Then, 
        \begin{align}
          \label{eq.theta.below}
          \lambda_a^+\int_{-\infty}^{\varepsilon_a/a} dy F_a^*(y) e^{\lambda_a^+y} dy &\ge a c_a + o(a c_a)\\
          \lambda_a^-\int_{-\infty}^{1/a} dy F_a^*(y) e^{\lambda_a^-y} dy &\le -a c_a + o(a c_a)
          \label{eq.theta.above}  
        \end{align} 
      and 
      \begin{align}
        \label{eq.theta.1moment.above}
        \int_{-\infty}^{1/a} dy  F_a^*(y) y e^{\lambda_a^+ y}\le \frac{\sigma^2}{2}+o(1),\quad a\to 0,\\
        \label{eq.theta.1moment.below}
        \int_{-\infty}^{\varepsilon_a/a} dy  F_a^*(y) y e^{\lambda_a^- y}\ge \frac{\sigma^2}{2}+o(1),\quad a\to 0.
      \end{align}  
    \begin{comment}
        Let $c_a=o(\theta_a)$. Then, 
        \begin{align}  
            \label{eq.theta.below}
            \e\left[e^{(\theta_a+c_a) X^{(a)}};X^{(a)}\le 1/a\right] \ge     1 + a c_a + o(a c_a), \quad a\to 0,\\ 
            \e\left[e^{(\theta_a-c_a) X^{(a)}};X^{(a)}\le 1/a\right] \le     1 - a c_a + o(a c_a), \quad a\to 0
            \label{eq.theta.above}.
        \end{align}   
      \end{comment}   
    \end{lemma}    
    \begin{proof}
    First, integrating by parts and using \eqref{eq.theta.below2}, 
      \begin{align*}
        &\lambda_a^+\int_{-\infty}^{1/a} dy F_a^*(y)e^{\lambda_a^+ y} =
        \lambda_a^+\left(
          -\int_{-\infty}^0 dy F_a(y) e^{\lambda_a^+ y} 
          + \int_{0}^{1/a} dy \overline F_a(y) e^{\lambda_a^+ y} 
        \right)\\ 
        &\hspace{1cm}=
          \overline F_a(1/a) e^{\lambda_a^+/a}-1 
          +\int_{-\infty}^0 dF_a(y) e^{\lambda_a^+ y} 
          + \int_{0}^{1/a} dF_a(y) e^{\lambda_a^+ y} 
        \\ 
        &\hspace{1cm}
        \ge \overline F_a(1/a) e^{\lambda_a^+/a}
        +ac_a+o(ac_a)\ge ac_a+o(ac_a).
    \end{align*}  
     Using the latter  inequality 
     and the assumption $\overline F(\varepsilon_a/a)=o(ac_a)$ 
     we obtain 
     \begin{align*}
      \lambda_a^+\int_{-\infty}^{\varepsilon_a/a} dy F_a^*(y)e^{\lambda_a^+ y}
      &\ge \lambda_a^+\int_{-\infty}^{1/a} dy F_a^*(y)e^{\lambda_a^+ y}
        -\lambda^+_a e^{\lambda^+/a} \overline F(\varepsilon_a/a)\\
      &\ge ac_a  +o(ac_a),
     \end{align*} 
     proving \eqref{eq.theta.below}. 

     Second, note that 
     for some constant $C=C(\sigma)$ depending only on $\sigma$, 
     \begin{align*}
      \int_{-\infty}^{1/a} dy  F_a^*(y) y e^{\lambda_a^+ y} 
      &=-\int_{-\infty}^0 y e^{\lambda_a^+ y} F_a(y)dy +\int_0^{1/a} y e^{\lambda_a^+ y} \overline F_a(y)dy \\ 
      &\le - \int_{-\infty}^0 yF_a(y)dy+\int_0^{1/a} y (1+C\lambda_a^+y)\overline F_a(y)dy \\ 
      &=\int_{-\infty}^\infty yF^*_a(y)dy - \int_{1/a}^\infty y\overline F_a(y)dy +C\lambda_a^+ \int_0^{1/a}y^2\overline F_a(y)dy\\ 
      &\le \frac{1}{2}\e[(X^{(a)})^2] +C\lambda_a^+ \int_0^{1/a}y^2\overline F(y)dy. 
     \end{align*} 
     It follows from  $\e[X^2]<\infty$ 
     that $y^2\overline F(y)\to 0$, 
     as  $y\to\infty$,  and hence 
     $$
     \lambda_a^+ \int_0^{1/a}y^2\overline F(y)dy = o(1).
     $$
     Since 
     $$
     \e[(X^{(a)})^2] = \sigma^2+a^2 
     =\sigma^2+o(1), 
     $$
     we obtain  \eqref{eq.theta.1moment.above}. 

     Similarly, integrating by parts 
     and using \eqref{eq.theta.above2} 
     \begin{align*}
      &\lambda_a^-\int_{-\infty}^{1/a} dy F_a^*(y)e^{\lambda_a^- y} 
      =
        \overline F_a(1/a) e^{\lambda_a^-/a}-1 
        +\int_{-\infty}^0 dF_a(y) e^{\lambda_a^- y} 
        + \int_{0}^{1/a} dF_a(y) e^{\lambda_a^- y} 
      \\ 
      &\hspace{1cm}
      \le \overline F_a(1/a) e^{\lambda_a^+/a}
      -ac_a+o(ac_a)\le -ac_a+o(ac_a)
  \end{align*}  
  we prove \eqref{eq.theta.above}. 

  Next, using the lower bound $e^x\ge 1+x,$ 
     \begin{align*}
      \int_{-\infty}^{1/a} dy  F_a^*(y) y e^{\lambda_a^- y} 
      &=-\int_{-\infty}^0 y e^{\lambda_a^- y} F_a(y)dy +\int_0^{1/a} y e^{\lambda_a^- y} \overline F_a(y)dy \\ 
      &\ge - \int_{-\infty}^0 y(1+\lambda^-_a y)F_a(y)dy+\int_0^{1/a} y \overline F_a(y)dy \\ 
      &=\int_{-\infty}^\infty yF^*_a(y)dy 
      - \int_{1/a}^\infty y\overline F_a(y)dy
      -\lambda_a^- \int_{-1/a}^0 y^2  F(y)dy
      \\ 
      &\ge \frac{1}{2}\e[(X^{(a)})^2] 
      - \int_{1/a}^\infty y\overline F_a(y)dy
      -\lambda_a^- \int_{-1/a}^0 y^2  F(y)dy
    \end{align*} 
     Since $\e[X^2]<\infty$ 
     we have the  convergence $y^2 F(y)\to 0$, as  $y\to-\infty$ and hence 
     $$
     \lambda_a^- \int_{-1/a}^0 y^2  F(y)dy= o(1).
     $$
     Next, 
     \begin{align*}
      \int_{-\infty}^{\varepsilon_a/a} dy  F_a^*(y) y e^{\lambda_a^- y} 
      &\ge \int_{-\infty}^{1/a} dy  F_a^*(y) y e^{\lambda_a^- y} 
      - e^{\varepsilon_a\lambda_a^- /a}\overline F(\varepsilon_a/a)\\ 
      &\ge \frac{\sigma^2}{2}-o(1), 
     \end{align*} 
    which proves \eqref{eq.theta.1moment.below}. 

    \end{proof}

    %%%%%%%%%%%%%%%%%%%%%%%%%%%%%%%%%%%%%%%%%%%%%%%%%%%%%%%%%%%%%%%%%%%%%%%%%%%%%%%%%%%%%%%%%%%%%%%%%%%%%%%%%%%%%%%%%%%%%%%

    \subsection{Decomposition of $\overline G_+$ and $\overline G_-$}
    To show sub/super-harmonic property 
    of functions $G_+$ and $G_-$
    we need to show  that 
    the differences 
    $$
    \beta_\pm(t) := \e \left[\overline G_\pm(t-X^{(a)})\right] - \overline G_\pm(t)
    $$
    are correspondingly positive/negative. 
    In this subsection 
    we will obtain a more convenient expression 
    for $G_\pm$ using the integration parts and then decompose 
    this expression. 
    To simplify notation we assume that $\sigma=1$.
    
    Both functions $\overline G_+(x)$ and 
    $\overline G_-(x)$ have a jump at $0$ 
    and are absolutely continuous on $[0,+\infty)$. 
    %For $x\neq 0$ let $g(x) = -\overline G'(x)$ and let $g(0)=0$.  
    For $x\ge 0$ let 
    \begin{equation}
      \label{def_G_density}
      g_\pm (x) := e^{\pm \alpha}\lambda^\pm_a e^{-\lambda^\pm_a x } 
      +2\lambda^\pm_a
      \int_0^x dz e^{-\lambda^\pm_a (x-z) } (x-z) \overline F_+ (z). 
    \end{equation}
    Then, it is not difficult to see that 
    $\overline G_\pm(x) = \int_x^\infty g_\pm(y) dy$ for $x\ge 0$.    
    Next, 
    \begin{align*}
      \beta_\pm(t) = e^{\pm 2\alpha} \overline F_a(t) 
      +\int_{-\infty}^t F_a(dy) \overline G_\pm (t-y) -  \overline G_\pm(t). 
    \end{align*}  
    The latter expression can be transformed using the integration by parts as follows, 
    \begin{align*}
      \beta_\pm(t) &= e^{\pm 2\alpha} \overline F_a(t) 
      +\int_{-\infty}^0 F_a(dy) \overline G_\pm (t-y)
      +\int_{0}^t F_a(dy) \overline G_\pm (t-y)-  \overline G_\pm(t)\\ 
      &= e^{\pm 2\alpha} \overline F_a(t)  
      +F_a(0)\overline G_\pm(t) +\overline F_a(0)\overline G_\pm(t)-\overline F_a(t)\overline G_\pm(0)
      -  \overline G_\pm(t)
      \\ 
      &+\int_0^t dy \overline  F_a(y) g_\pm(t-y)  -\int_{-\infty}^0 dy F_a(y) g_\pm(t-y). 
    \end{align*}  
    Using \eqref{def_f_star} 
    we can write a more compact expression for $\beta_\pm $ as follows, 
    \begin{equation}
      \label{beta_int_by_parts}
      \beta_\pm(t) 
      = (e^{\pm 2\alpha}-\overline G_\pm(0)) \overline F_a(t) + 
      \int_{-\infty}^t dy F_a^*(y) g_\pm(t-y).
    \end{equation}  
    It follows from \eqref{def_G} that 
    $$
    \overline G_\pm(0)= e^{\pm \alpha} +\frac{2}{\lambda^\pm_a\sigma^2} \overline F^I(\varepsilon_a/a). 
    $$
    Since $\alpha$ is fixed and by assumptions of Theorem~\ref{thm.main} 
    $\overline F^I(\varepsilon_a/a)=o(a)$,  
    for all  sufficiently small $a$, 
    \begin{equation}
      \label{g_at_zero}
    \overline G_+(0) \le \frac{1}{2}e^{2\alpha}+ \frac{1}{2}e^{\alpha} , 
    \quad  
    \overline G_-(0) \ge  \frac{1}{2}e^{-\alpha}+\frac{1}{2}e^{-2\alpha}.  
    \end{equation} 
    In what follows we will use the following 
    decomposition 
    $$
    \beta_\pm(t) = I^\pm (t)+E^\pm(t)+
    (e^{\pm \alpha}-\overline G_\pm(0)) \overline F_a(t), 
    %(\overline G_\pm(0)-1) \overline F(t),
    $$
    where 
    \begin{align*}
      I^\pm(t)&:= 
      \lambda_a^\pm e^{\pm \alpha}\int_{-\infty}^t dy F_a^*(y) 
      e^{-\lambda_a^\pm (t-y)}\\ 
          &=
          \lambda_a^\pm e^{\pm \alpha}\left(\int_{-\infty}^{\varepsilon_a/a} 
            +\int_{\varepsilon_a/a}^{t} 
            \right)
            dy F_a^*(y) e^{-\lambda_a^\pm (t-y)}\\
             &=: I^\pm_1(t)+I^\pm_2(t)
    \end{align*}  
    and 
    \begin{align*}
      E^\pm(t)&:= 2\lambda_a^\pm\int_{-\infty}^{\varepsilon_a/a} dy  F_a^*(y) \int_0^{t} dz e^{-\lambda_a^\pm(t-y-z)}(t-y-z)\overline F_+(z)\\
          &\hspace{1cm}-2\lambda_a^\pm\int_{0}^{\varepsilon_a/a} dy  F_a^*(y) \int_{t-y}^{t} dz e^{-\lambda_a^\pm(t-y-z)}(t-y-z)\overline F_+(z)\\  
          &\hspace{1cm}+2\lambda_a^\pm\int_{-\infty}^{0} dy  F_a^*(y) \int_{t}^{t-y} dz e^{-\lambda_a^\pm(t-y-z)}(t-y-z)\overline F_+(z)\\
          &\hspace{1cm}+2\lambda_a^\pm\int_{\varepsilon_a/a}^t dy  \overline F_a(y) \int_{0}^{t-y} dz e^{-\lambda_a^\pm(t-y-z)}(t-y-z)\overline F_+(z)\\
          &=:E^\pm_1(t)+E^\pm_2(t)+E^\pm_3(t)+E^\pm_4(t).
    \end{align*}

    \subsection{Subharmonic property of $\overline G_+$}
    In this subsection we will 
    prove the subharmonic property of $\overline G_+$,  
    that is that  $\beta_+(t)\ge 0$,   for all  $t\ge \varepsilon_a/a$ and 
    for all sufficiently small $a$.  
    \begin{lemma}
      \label{lem.i1.i2.e1}
      Assume that $\overline F(\varepsilon_a/a)=o(ac_a)$. 
      Then, there exists $a_0$ such that for $t\ge \varepsilon_a/a$ and $a\in(0,a_0)$
      \begin{equation}
        \label{i1.i2.e1}
        I^+(t)+E_1^+(t) \ge 0. 
      \end{equation}  
    \end{lemma}  
    \begin{proof}
    First we  apply \eqref{eq.theta.below} to obtain 
    uniformly in $t\ge 0$, 
    \begin{equation}
      \label{i1}
      I_1^+(t) = e^{-\lambda_a^+ t}\lambda_a^+e^\alpha 
      \int_{-\infty}^{\varepsilon_a/a} dy F_a^*(y) e^{\lambda_a^+ y }
      \ge e^\alpha e^{-\lambda_a^+ t}(ac_a+o(ac_a)).
      \end{equation}  
     Next, changing the order of integration 
     and using \eqref{eq.theta.below} and \eqref{eq.theta.1moment.above} 
     we obtain 
     \begin{align*}
      E_1^+(t)&=
      2\lambda_a^+
      \int_0^{t} dze^{-\lambda_a^+(t-z)}\overline F_+(z)
      \int_{-\infty}^{\varepsilon_a/a} dy  F_a^*(y)  e^{\lambda_a^+ y}(t-y-z)\\
      &\ge 
      \int_0^{t} dze^{-\lambda_a^+(t-z)}\overline F_+(z) 
      \left(
      %2(t-z)ac_a  
      -2\lambda_a^+
      \left(\frac{1}{2}-o(1)\right)\right)\\
      &\ge (-\lambda_a^+-o(\lambda_a^+))\int_0^{t} dze^{-\lambda_a^+(t-z)}\overline F_+(z). 
     \end{align*}    
     Therefore, for a sufficiently small $a_1$ and $a<a_1$ we obtain 
     \begin{align*}
      I_2^+(t)+E_1^+(t) 
      &\ge 
      ((e^\alpha-1)\lambda_a^+-o(\lambda_a^+))\int_{\varepsilon_a/a}^{t} dze^{-\lambda_a^+(t-z)}\overline F_+(z)
      \\
      &\hspace{1cm}+
      (-\lambda_a^+-o(\lambda_a^+))\int_0^{\varepsilon_a/a} dze^{-\lambda_a^+(t-z)}\overline F_+(z)\\
      &\ge
      (-\lambda_a^+-o(\lambda_a^+))\int_0^{\varepsilon_a/a} dze^{-\lambda_a^+(t-z)}\overline F_+(z)\\
      &\ge 
      (-1-o(1)) \overline F(\varepsilon_a/a)e^{-\lambda_a^+t}. 
      \end{align*} 
      Since $\overline F(\varepsilon_a/a) = o(ac_a)$ 
      the inequality \eqref{i1.i2.e1}
      follows 
      from the latter inequality and  
      \eqref{i1}. 
    \end{proof}
    \begin{lemma}
      \label{lem.sub}
      Assume that $\overline F(\varepsilon_a/a)=o(ac_a)$. 
      Then, there exists $a_0$ such that for 
      $t\ge \varepsilon_a/a$ and $a\in(0,a_0)$, 
      \begin{equation}
        \label{beta.plus}
        \beta_+(t)\ge 0. 
      \end{equation}  
    \end{lemma} 
    \begin{proof}
     Clearly $E_2^+(t)\ge 0$ and $E_4^+(t)\ge 0$. % and $I_4(t)\ge 0$. 
     Then, in view of  Lemma~\ref{lem.i1.i2.e1},  it is sufficient to show that 
     uniformly in $t\ge \varepsilon_a/a$ for all sufficiently small $a$ 
     \begin{equation}
        \label{e3}
     E_3^+(t)+(e^{2\alpha}-\overline G_+(0))\overline F(t) \ge 0. 
     \end{equation}
    First changing the order of integration and then the variables  we obtain 
     \begin{align*}
        E_3^+(t) &= -2\lambda_a^+\int_{-\infty}^{0} dy   F_a(y) \int_{t}^{t-y} dz e^{-\lambda_a^+(t-y-z)}(t-y-z)\overline F_+(z)\\ 
               &= -2\lambda_a^+ \int_{t}^{\infty} dz \overline F_+(z)
               \int_{-\infty}^{t-z} 
               dy e^{-\lambda_a^+(t-y-z)}(t-y-z) F_a(y)
               \\ 
               &=-2\lambda_a^+ \int_{t}^{\infty} dz \overline F_+(z)
               \int_{0}^{\infty}dy e^{-\lambda^+_a y} y F_a(t-z-y)dy. 
     \end{align*} 
     Splitting the outer integral in two parts we obtain 
     \begin{align*}
        E_3^+(t)&= -2\lambda_a^+ 
        \left(\int_{t}^{t+\varepsilon_a/a} 
        +\int_{t+\varepsilon_a/a}^\infty 
        \right)
        dz \overline F_+(z)
        \int_{0}^{\infty}dy e^{-\lambda^+_a y} y F_a(t-z-y)dy\\ 
        &\ge \frac{-2\lambda^+_a \varepsilon_a}{a}
              \overline F_+(t)  \int_0^\infty yF_a(-y)dy 
              -2\lambda_a^+ \overline F_+(t) 
              \int_{t+\varepsilon_a/a}^\infty F_a(t-z) dz 
              \int_{0}^{\infty}dy e^{-\lambda^+_a y} y\\ 
         &\ge      
         \overline F_+(t) 
         \left(
            -\frac{2\lambda^+_a \varepsilon_a}{a} 
            \int_0^\infty yF_a(-y)dy
            -2\int_{-\infty} ^{-\varepsilon_a/a} F_a(z)dz 
         \right)
         \ge o(1)\overline F_+(t).
     \end{align*}
     In view of \eqref{g_at_zero} and since 
     $\alpha$ is fixed, 
     we  obtain \eqref{e3} for all sufficiently small $a$.    
\end{proof}

    \subsection{Superharmonic property of $\overline  G_-$}
    In this subsection we will 
    prove the superharmonic property of $G_-$, 
    that is that there exists $a_0$ such that 
    $\beta_-(t)\le 0$ for all $a\in (0,a_0)$ 
    and $t\ge \varepsilon_a/a$.  
    
    \begin{lemma}
      \label{lem.i1.i2.e1.super}
      Assume that $\overline F(\varepsilon_a/a)=o(ac_a)$ 
      and condition \eqref{subex} holds. 
      Then, there exists $a_0$ such that 
      for $t\ge \varepsilon_a/a$ and $a\in(0,a_0)$
      \begin{equation}
        \label{i1.i2.e1.super}
        I^-(t)+E_1^-(t) +E_4^-(t) \le 0.%o(c_a)\overline F(t). 
      \end{equation}  
    \end{lemma}  
    \begin{proof}
    First we  apply \eqref{eq.theta.above} to obtain 
    uniformly in $t\ge 0$, 
    \begin{equation}
      \label{i1.super}
      I_1^-(t) = e^{-\lambda_a^- t}\lambda_a^-e^{-\alpha} \int_{-\infty}^{\varepsilon_a/a} dy F_a^*(y) e^{\lambda_a^- y }
      \le e^{-\alpha} e^{-\lambda_a^- t}(-ac_a+o(ac_a))
      \le 0.
      \end{equation}  
     Next, changing the order of integration 
     and using \eqref{eq.theta.above} and \eqref{eq.theta.1moment.below} 
     we obtain 
     \begin{align*}
      E_1^-(t)&=
      2\lambda_a^-
      \int_0^{t} dze^{-\lambda_a^-(t-z)}\overline F_+(z)
      \int_{-\infty}^{\varepsilon_a/a} dy  F_a^*(y)  e^{\lambda_a^- y}(t-y-z)\\
      &\le 
      \int_0^{t} dze^{-\lambda_a^-(t-z)}\overline F_+(z) 
      \left(
      -2(t-z)ac_a  
      -2\lambda_a^-
      \left(\frac{1}{2}+o(1)\right)\right)\\
      &\le (-\lambda_a^-+o(\lambda^-_a))\int_0^{t} dze^{-\lambda_a^-(t-z)}\overline F_+(z). 
     \end{align*}    
     Also, changing the variable and then the order of integration we obtain 
     \begin{align*}
      E_4^-(t)&=2\lambda_a^-\int_{\varepsilon_a/a}^t dy  \overline F_a(y) \int_{0}^{t-y} dz e^{-\lambda_a^-z}z\overline F_+(t-y-z)\\
      &=2\lambda_a^-\int_0^{t-\varepsilon_a/a}dz  e^{-\lambda_a^-z}z \int_{\varepsilon_a/a}^{t-z} dy \overline F_a(y)\overline F_+(t-y-z) 
     \end{align*} 
     Using the assumption \eqref{subex}
     we obtain that 
    \begin{align*}
      E_4^-(t)&\le 
      o(c_a)2\lambda_a^-\int_0^{t-\varepsilon_a/a}dz  e^{-\lambda_a^-z}z\overline F_+(t-z)\\ 
            &=o(ac_a)2\int_{\varepsilon_a/a}^{t}dz  e^{-\lambda_a^-(t-z)}(t-z)\overline F_+(z)
    \end{align*}  
     Therefore, for a sufficiently small $a_0$ and $a<a_0$ we obtain 
     \begin{align*}
      I_2^-(t)+E_1^-(t) + E_4^-(t)
      &\le 
      ((e^{-\alpha}-1)\lambda^-_a+o(\lambda_a^-))
      \int_{\varepsilon_a/a}^{t} dze^{-\lambda_a^-(t-z)}\overline F_-(z)
      \\
      %&\hspace{1cm}+
      %2\lambda_a^-o(ac_a)\int_{t-\varepsilon/a}^{t} 
      %dze^{-\lambda_a^-(t-z)} (t-z)\overline F_+(z)\\
      %&\le
      %2\lambda_a^-o(ac_a)\int_{t-\varepsilon/a}^{t} 
      %dze^{-\lambda_a^-(t-z)} (t-z)\overline F_+(z)
      %\\
      %%&\le 
      %2\lambda_a^-o(ac_a) \left(\frac{\varepsilon}{a}\right)^2 
      %\overline F_+(t-\varepsilon/a)\\ 
      &\le 0.%  o(c_a)\overline F(t).
    \end{align*} 
      %Since $\overline F(\varepsilon/a) = o(ac_a)$ 
      Then  the inequality \eqref{i1.i2.e1.super}
      follows 
      from the latter inequality and  
      \eqref{i1.super}. 
    \end{proof}

    \begin{lemma}
      \label{lem.super}
      Assume that $\overline F(\varepsilon_a/a)=o(ac_a)$. 
      Then, there exists $a_0$ such that for 
      $t\ge  \varepsilon_a/a$ and $a\in(0,a_0)$
      \begin{equation}
        \label{beta.minus}
        \beta_-(t)\le 0. 
      \end{equation}  
    \end{lemma} 
    \begin{proof}
     Since $E_3^-(t)\le 0$,
     in view of    Lemma~\ref{lem.i1.i2.e1.super},  
     it is sufficient to show that 
     uniformly in $t\ge \varepsilon_a/a$ 
     for all sufficiently small $a$ 
     \begin{equation}
        \label{e2}
     E_2^-(t)+(e^{-2\alpha}-\overline G_-(0))\overline F(t) \le 0. 
     \end{equation}
    Changing the order of integration we obtain 
     \begin{align*}
        E_2^-(t)& = 
        -2\lambda_a^-
          \int_{t-\varepsilon_a/a}^t
          dz e^{-\lambda_a^-(t-z)} \overline F_+(z)
        \int_{t-z}^{\varepsilon_a/a} e^{\lambda_a^- y} F_a^*(y)(t-y-z)dy\\
    &=
      2\lambda_a^-\int_{t-\varepsilon_a/a}^t 
        dz e^{-\lambda_a^-(t-z)} \overline F_+(z)
        \int_{t-z}^{\varepsilon_a/a} e^{\lambda_a^- y} F_a^*(y)(y-(t-z))dy\\ 
        &\le 4\lambda_a^- e^{\varepsilon_a \lambda_a^-/a } 
        \left(
        \frac{\varepsilon_a}{a}
        \right)
        \int_{t-\varepsilon/a}^t 
        dz  \overline F_+(z)
        \int_{t-z}^{\varepsilon_a/a}  \overline F_a(y)dy\\ 
        &\le 4\lambda_a^- e^{\varepsilon \lambda_a^-/a }  \left(
          \frac{\varepsilon_a}{a}
          \right)
          \overline F(t)
        \int_0^{\varepsilon_a/a} du \int_u^{\varepsilon_a/a} \overline F_a(y)dy
        . 
     \end{align*} 
     Since  $\varepsilon_a\to 0$ 
     and $\alpha$ is fixed using \eqref{g_at_zero} 
     we obtain \eqref{e2}. 
    \end{proof}

    \subsection{Proof of Theorem~\ref{prop1}}
    Let 
    $$
    \mu^{(a)}(x):=\min\{n\ge 1: S_n^{(a)}>x\}.
    $$
    Then, clearly  
    $
    \{M^{(a)}>x\} = \{\mu^{(a)}(x)<\infty\}.
    $
    \begin{proof} (of Theorem~\ref{prop1})
      We start with the lower bound. 
      Consider the function 
      $
      \widehat G_+(t) = \overline G_+(t+\varepsilon_a). 
      $
      By Lemma~\ref{lem.sub} for $x\ge 0$ 
      and all sufficiently small $a$, 
      $$
      \e[\widehat G_+(x-X^{(a)})]-\widehat G_+(x) \ge 0. 
      $$
      Consequently, 
      $$
      \widehat G_+(x-S_{n\wedge \mu^{(a)}(x)}^{(a)}) 
      \mbox { is a bounded non-negative submartingale.}
      $$
      Then, by the optional stopping theorem, 
      \begin{align*}
        \widehat G_+(x) 
        &\le \e[\widehat G_+(x-S_{\mu^{(a)}(x)}^{(a)}) ] 
        = \e[\widehat G_+(x-S_{\mu^{(a)}(x)}^{(a)});\mu^{(a)}(x)<\infty] \\ 
        &\le e^{2\alpha} \pr(\mu^{(a)}(x)<\infty)
        =e^{2\alpha} \pr(M^{(a)}>x).
      \end{align*}  
      This proves the lower bound \eqref{eq.lower.bound}. 

      The proof of the upper bound is analogous. 
      Consider the function 
      $
      \widehat G_-(t) = \overline G_-(t+\varepsilon_a). 
      $
      By Lemma~\ref{lem.super} for $x\ge 0$ 
      and all sufficiently small $a$, 
      $$
      \e[\widehat G_-(x-X^{(a)})]-\widehat G_-(x) \le 0. 
      $$
      Consequently, 
      $$
      \widehat G_-(x-S_{n\wedge \mu^{(a)}(x)}^{(a)}) 
      \mbox { is a bounded non-negative supermartingale.}
      $$
      Then, by the optional stopping theorem, 
      \begin{align*}
        \widehat G_-(x) 
        &\ge \e[\widehat G_-(x-S_{\mu^{(a)}(x)}^{(a)}) ] 
        = \e[\widehat G_-(x-S_{\mu^{(a)}(x)}^{(a)});\mu^{(a)}(x)<\infty] \\ 
        &\ge e^{-2\alpha} \pr(\mu^{(a)}(x)<\infty)
        =e^{-2\alpha} \pr(M^{(a)}>x).
      \end{align*}  
      This proves the lower bound \eqref{eq.upper.bound}.

    \end{proof}

    \section{Proof of theorem~\ref{thm.main}}
    \label{sec.main}
    We will derive Theorem~\ref{thm.main} from Theorem~\ref{prop1}. 
    For that we need first to define suitable $c_a$ and $\varepsilon_a$, 
    which will be done in Lemma~\ref{lem.sequences} below. 
    Then, we will show that the error terms $c_a$ and $\varepsilon_a/a$ 
    are negligible and can be removed from the asymptotics. 
    This will  be done via a sequence of Lemmas. 
    To simplify notation we will assume throughout that $\sigma=1$.
  
    We will repeatedly use the following 
    bound that  follows from \eqref{sc2} and  \eqref{sc3} 
    and is valid  for any $z>0$ 
    and all sufficiently large $x$, 
    \begin{align}
      \label{sc3-cor}
      g(x+z)-g(x)=\int_x^{x+z}g'(u)du
      \le \int_x^{x+z} \gamma_0 \frac{g(u)}{u}du 
      \le \gamma_0 g(x)\frac{z}{x}.
    \end{align}  
    Also, it follows from 
    \eqref{eq.boundary} 
    and \eqref{sc2} 
    that 
    \begin{align}
      \label{eq.bound1.xa1}
      a x(a)  &= \frac{\sigma^2}{2} g(x(a)) \uparrow
       \infty,\quad a\to 0.
       %\\  
      %\label{eq.bound1.xa2}
      %a x(a)^{1-\gamma_0}  
      %&=\frac{\sigma^2}{2} \frac{g(x(a))}{x(a)^{\gamma_0}} 
      %\downarrow 0. 
    \end{align}   

    \begin{lemma}
      \label{lem.sequences}
      Let $F$ satisfy \eqref{sc1} and \eqref{sc2} 
      and $x(a)$ solve \eqref{eq.boundary}. 
      Then there exist functions  $\varepsilon_a$ and $c_a$ 
      such that $c_a=o(1/x(a))$, $\varepsilon_a\to 0$, 
      $\varepsilon_a/a\to +\infty$, 
      $
      \overline F(\varepsilon_a/a) = o(ac_a)
      $
      and \eqref{subex} holds. 
  \end{lemma}  
  \begin{proof}
    It follows from \eqref{eq.bound1.xa1} that 
    $x(a)\ge \varepsilon/a$  for any fixed $\varepsilon$, as $a\to 0$.  
    Now   observe  
       that \eqref{eq.boundary} 
       and \eqref{sc2} imply that 
       for any $\varepsilon>0$
       \begin{align*}
        2a &= \frac{g(x(a))}{x(a)} 
        =
        \frac{g(x(a))}{x(a)^{\gamma_0}} x(a)^{\gamma_0-1}
        \le 
        \frac{g(\varepsilon/a)}{(\varepsilon/a)^{\gamma_0}} x(a)^{\gamma_0-1}
       \end{align*} 
       and, consequently, 
       \begin{equation}
        \label{eq.domin.last.3}
        ax(a) \varepsilon^{\frac{\gamma_0}{1-\gamma_0}}
       \le
       \left(\frac{g(\varepsilon/a)}{2}\right)^{\frac{1}{1-\gamma_0}}. 
       \end{equation}

       Now  let $\delta_0>0$ be such that $\gamma_0(1+\delta_0)<1$ and 
       let 
       \begin{equation}
        \label{def.varepsilon.a}
        \varepsilon_a = (ax(a))^{-\delta_0},
       \end{equation} 
       which is  monotone decreasing to $0$, as $a\to 0$,  
       by \eqref{eq.bound1.xa1}. 
       Moreover, by \eqref{eq.boundary}, 
       \begin{align*}
       \frac{\varepsilon_a}{a} &= 
       a^{-\delta_0-1}x(a)^{-\delta_0}
       =
       x(a) \left(\frac{\sigma^2}{2} g(x(a))\right)^{-\delta_0-1}
       \\
       &= x(a)^{1-\gamma_0(1+\delta_0)} 
       \left(\frac{\sigma^2}{2} \frac{g(x(a))}{x(a)^{\gamma_0}}\right)^{-\delta_0-1}
       \uparrow \infty, 
       \end{align*}
       as $a\to 0$, by \eqref{sc2}. 
       Also, 
       plugging in $\varepsilon=\varepsilon_a$  in 
       \eqref{eq.domin.last.3} 
       we obtain
       \begin{equation}
        \label{axa}
       (ax(a))^{1-\gamma_0(1+\delta_0)} \le 
       \frac{g(\varepsilon/a)}{2}. 
       \end{equation}
      Next, uniformly in $x\ge 2\varepsilon_a/a$ 
       \begin{align*}
        I&:=\int_{\varepsilon_a/a}^{x-\varepsilon_a/a}
        \frac{\overline F(y)\overline F(x-y)}{\overline F(x)}dy 
        = 
        2\int_{\varepsilon_a/a}^{x/2}
        \frac{\overline F(y)\overline F(x-y)}{\overline F(x)}dy 
        \\
        &\sim 
        2\int_{\varepsilon_a/a}^{x/2}
        \frac{x^2}{y^2(x-y)^2} e^{g(x)-g(x-y)-g(y)}dy 
        \le 8 
        \int_{\varepsilon_a/a}^{x/2}
        e^{g(x)-g(x-y)-g(y)}y^{-2}dy. 
      \end{align*} 
      Since $y\le x/2$, 
      applying \eqref{sc3-cor} and then \eqref{sc2}
      we obtain 
      \begin{equation}
        \label{subex-diff}
      g(x)-g(x-y)-g(y)
      \le \gamma_0 y \frac{g(x-y)}{x-y} - g(y)
      \le (\gamma_0-1) g(y). 
      \end{equation}
      Therefore, 
      \begin{align*}
        I\le 
        8 
        \int_{\varepsilon_a/a}^{x/2}
        e^{(\gamma_0-1) g(y)} y^{-2}dy 
        \le 
        8 e^{(\gamma_0-1) g(\varepsilon_a/a)} 
        \int_{\varepsilon_a/a}^{x/2}
        y^{-2}dy 
        \le 8 \frac{a}{\varepsilon_a} e^{(\gamma_0-1) g(\varepsilon_a/a)}. 
      \end{align*}  
      Next, we will show that $I=o(1/x(a))$ uniformly in $x\ge 2\varepsilon_a/a$. 
      For that note 
      using the definition of $\varepsilon_a$ and 
      \eqref{axa}
      \begin{align*}
      \frac{I}{(1/x(a))}
      &\le 8 \frac{a x(a)}{\varepsilon_a} e^{(\gamma_0-1) g(\varepsilon_a/a)}
      = 8 (a x(a))^{1+\delta_0}e^{(\gamma_0-1) g(\varepsilon_a/a)}\\
      &\le 
      \left(
        \frac{g(\varepsilon/a)}{2}
      \right)^{
        \frac{1+\delta_0}{1-\gamma_0(1+\delta_0)}
      }e^{(\gamma_0-1) g(\varepsilon_a/a)} \to 0, 
      \end{align*}
      since $g(\varepsilon_a/a)\to \infty $ as $a\to 0 .$
      Since $I=o(1/x(a))$ there exists 
      a function $\widetilde c_a \uparrow \infty$ such that 
      $\widetilde c_a=o(1/x(a))$ and still $I=o(\widetilde c_a). $
      
      Now note that 
      \begin{align*}
        \frac{\overline F(\varepsilon_a/a)}
        {a/x(a)}
        &\sim \frac{a x(a)}{\varepsilon_a^2}
        e^{-g(\varepsilon_a/a)}
        \le 
        (a x(a))^{1+2\delta_0}e^{-g(\varepsilon_a/a)}\\
        &\le \left(
        \frac{g(\varepsilon/a)}{2}
      \right)^{
        \frac{1+2\delta_0}{1-\gamma_0(1+\delta_0)}
      }e^{- g(\varepsilon_a/a)}\to 0, 
       \quad a\to 0, 
      \end{align*}  
      since $g(\varepsilon_a/a)\to \infty $ as $a\to 0 .$
     Since $\overline F(\varepsilon_a/a) =o(a/x(a)) $ 
     there exists $\widehat c_a$ such that 
     $\widehat c_a=o(1/x(a))$  and still $\overline F(\varepsilon_a/a)=o(a\widehat c_a), $
      as $a\to 0$. Now we can simply put $c_a=\max(\widetilde c_a,\widehat c_a).$ 
  \end{proof}

    We will need the following insensitivity property. 
    \begin{lemma}(Insensitivity)
      \label{lem.insensitivity}
      Let $F$ satisfy \eqref{sc1} and \eqref{sc2}. 
      Let $c_a = o(1/x(a))$%, \theta_a\sim \frac{2a}{\sigma^2}$ 
      and $\varepsilon_a\to 0$. Then,  for any fixed $A>1$, 
      as $a\to 0$, 
      \begin{align}
      \label{exp-long-tails}
      e^{-(\theta_a\pm c_a)x} &\sim e^{-\theta_a x},&& \mbox{    uniformly in  $x\le Ax(a)$, }\\ 
      \label{subexp-long-tails}
      \overline F\left(
         x+\frac{\varepsilon_a}{a}
        \right)&\sim 
      \overline F(x), 
      &
      \overline F^I\left(
         x+\frac{\varepsilon_a}{a}
        \right)\sim 
      \overline F^I(x), 
      & \mbox{    uniformly in  $x\ge x(a)/A$, }
      \end{align}  
    \end{lemma}  
    \begin{proof}
      Property \eqref{exp-long-tails} follows immediately from the fact 
      that uniformly in $x\le Ax(a)$ 
      $$
      1\le e^{c_ax} \le e^{Ac_ax(a)} \to 1, \quad a\to 0.
      $$
      Next, using \eqref{sc2} and \eqref{sc3-cor},
      \begin{align*}
        1&\le 
        \frac
        {\overline F\left(
          x
         \right)}
        {\overline F\left(
          x+\frac{\varepsilon_a}{a}
         \right)}
         \sim 
         \left(
           \frac{x+\varepsilon_a/a}{x}
         \right)^2 
         \exp\left\{g\left(x+\frac{\varepsilon_a}{a}\right)-g(x)\right\} \\ 
         &\sim 
         \exp\left\{g\left(x+\frac{\varepsilon_a}{a}\right)-g(x)\right\} \\ 
        &\le 
        \exp\left\{
        \gamma_0\frac{\varepsilon_a}{a} \frac{g(x)}{x}
        \right\}
        \le 
        \exp\left\{
        \gamma_0\frac{\varepsilon_a}{a} \frac{g(x(a)/A)}{x(a)/A}
      \right\}
      \\ 
     &
     \le 
      \exp\left\{
      \gamma_0A\frac{\varepsilon_a}{a} \frac{g(x(a))}{x(a)}
     \right\}
     = 
     \exp\left\{
      \gamma_0A\frac{\varepsilon_a}{a} \frac{2a}{\sigma^2}
     \right\}\to 1.
       \end{align*}  
      The second equivalence in \eqref{subexp-long-tails} 
      follows immediately from the first, since 
      uniformly in $x\le Ax(a)$,  
      $$
      \overline F^I\left(
         x+\frac{\varepsilon_a}{a}
        \right) = 
      \int_x^\infty
      F\left(y+\frac{\varepsilon_a}{a}\right)dy 
      \sim 
      \int_x^\infty
      F(y)dy  = \overline F^I(x).
      $$
    \end{proof}  
    We will also require the following upper and lower bounds for the 
    integrated tail.
    \begin{lemma}
      \label{lem.bounds.integrated.tail}
      Let $F$ satisfy \eqref{sc1} and \eqref{sc2}. Then, 
      for any $\delta>0$
      there exists $x_0$  such that 
      for $x>x_0,$
      \begin{align}
        \label{eq.bounds.integrated.tail.upper}
        \overline F^I(x) &\le (1+\delta)x^{-1}e^{-g(x)}\\ 
        \label{eq.bounds.integrated.tail.lower}
        \overline F^I(x)&
        %\ge 
        %c \gamma_0  (x^{2}g'(x))^{-1} e^{-g(x)}
        \ge \frac{1}{2}e^{-\gamma_0}(xg(x))^{-1}e^{-g(x)}.
      \end{align}  
    \end{lemma}   
    \begin{proof}
      For the upper bound  note that, as $x\to \infty$, 
      $$
      \overline F^I (x)   
      =\int_x^{\infty} \overline F(y)dy 
      \sim 
      \int_x^{\infty} y^{-2}e^{-g(y)}dy 
      \le  e^{-g(x)}
      \int_x^{\infty}y^{-2}dy
      =x^{-1}e^{-g(x)}. 
      $$
      For the lower bound  note that, as $x\to\infty$,  
     \begin{align}
      \nonumber
      \overline F^I(x) &\ge 
      \int_x^{x+x/g(x)}\overline F(y) dy 
      \sim\int_x^{x+x/g(x)}y^{-2}e^{-g(y)}dy \\ 
      \label{eq.lower.bound.1}
      &\ge e^{-g(x+x/g(x))} \int_x^{x+x/g(x)}y^{-2}dy
      \sim e^{-g(x+x/g(x))} \frac{1}{xg(x)}
      \end{align} 
      Applying \eqref{sc3-cor} we obtain, 
      \begin{align*}
        g(x+x/g(x)) - g(x)\le 
        \gamma_0.
      \end{align*}  
      Plugging in the latter inequality in 
       \eqref{eq.lower.bound.1}
      that 
      $$
      \overline F^I(x) 
      \ge (1+o(1)) 
      \frac{1}{xg(x)} e^{-g(x)} e^{-\gamma_0},
      $$
      which implies  \eqref{eq.bounds.integrated.tail.lower}.

    \end{proof}  
    Other useful bounds are given in the following Lemma.
    \begin{lemma}
      \label{lem.bounds.via.xa}
      Let $x(a)$ be solution to 
      \eqref{def.theta.xa}. Let $x_0$ be such that 
      $g(x)/x^{\gamma_0}$ is decreasing for $x>x_0$, where 
      $\gamma_0\in (0,1)$. Then, 
      \begin{align}
        \label{eq.bounds.via.xa.above}
        \frac{x}{g(x)} &\le \frac{\sigma^2}{2a} 
        \left(
          \frac{x}{x(a)}
        \right)^{1-\gamma_0},\quad x_0\le x\le x(a),\\ 
        \label{eq.bounds.via.xa.below}
        \frac{x}{g(x)} &\ge \frac{\sigma^2}{2a} 
        \left(
          \frac{x}{x(a)}
        \right)^{1-\gamma_0}, \quad  x\ge x(a).
      \end{align}  
    \end{lemma}  
    \begin{proof}
      We have, for $x: x_0\le x\le x(a)$, 
    \begin{align*}
        \frac{x}{g(x)} 
        =\frac{x^{\gamma_0}}{g(x)}x^{1-\gamma_0}
        \le  \frac{x(a)^{\gamma_0}}{g(x(a))}x^{1-\gamma_0}
         = \frac{x(a)}{g(x(a))}
         \left(\frac{x}{x(a)}\right)^{1-\gamma_0}
         =\frac{\sigma^2}{2a}
         \left(\frac{x}{x(a)}\right)^{1-\gamma_0}. 
      \end{align*}  
      Similarly, for $x\ge x(a)$, 
      \begin{align*}
        \frac{x}{g(x)} 
        =\frac{x^{\gamma_0}}{g(x)}x^{1-\gamma_0}
        \ge  \frac{x(a)^{\gamma_0}}{g(x(a))}x^{1-\gamma_0}
         =\frac{\sigma^2}{2a}
         \left(\frac{x}{x(a)}\right)^{1-\gamma_0}. 
      \end{align*}  
      
    \end{proof}  

    Next we will show that for $x\le (1-\delta)x(a)$ 
    the exponential term dominates over subexponential 
    in $\overline G_\pm$. 
    \begin{lemma}
      \label{lem.exp-domination}
      Let $\varepsilon_a\to 0$ so that 
      $\overline F^I(\varepsilon_a/a) = o(a)$
      and let $\theta_a\sim 2a/\sigma^2$ as $a\to 0$.
      Then, for any $\delta\in(0,1)$, 
      uniformly in $x\le (1-\delta)x(a)$, 
      \begin{equation}
        \label{exp-dominates}
        \frac{2}{\theta_a }
        \overline F^I(x\vee \varepsilon_a/a)
        =o(e^{-\theta_a x}). 
      \end{equation}  
    \end{lemma} 
    \begin{proof}
      To simplify notation we assume that $\sigma=1$. 
      Equivalence \eqref{exp-dominates} clearly holds 
      uniformly in  $x\le 1/\theta_a$, 
      since $\overline F^I(\varepsilon_a/a) = o(a)$.
      Hence it is sufficient to prove \eqref{exp-dominates} 
      for $x\ge 1/\theta_a$.
      Since $\lambda^\pm_a\sim \frac{2a}{\sigma^2}$, 
      by applying \eqref{eq.bounds.via.xa.above}, we obtain 
      for $x\in [1/\theta_a,(1-\delta)x(a)]$, 
      \begin{align}
        \nonumber 
        -\lambda^\pm_a x +g(x) 
        &=g(x)\left(1-
        \lambda^\pm_a\frac{x}{g(x)} 
        \right) 
        \ge 
        g(x)
        \left(
        1-\frac{\lambda^\pm_a}{2a} \left(\frac{x}{x(a)}\right)^{1-\gamma_0}
        \right)
        \\ 
        &
        \label{eq.minimizer}
        \ge g(x)
        \left(
          1-\frac{\lambda^\pm_a}{2a} (1-\delta)^{1-\gamma_0}
        \right)
        \ge \delta_1g(x),
      \end{align}  
      for some $\delta_1 \in (0,1-(1-\delta)^{1-\gamma_0}) $ 
      and all sufficiently small $a$. 
      
      Therefore, using  \eqref{eq.bounds.integrated.tail.upper} and \eqref{eq.minimizer} 
      we obtain, 
      as $a\to 0$, 
      \begin{align*}
        \frac{e^{-\theta_a x}}{\frac{2}{\theta_a}\overline F^I(x)}
        &\ge \frac{1}{4}\frac{\theta_a e^{-\theta_a x}}{x^{-1}e^{-g(x)}}
        =\frac{1}{4}\exp
        \left\{
          -\theta_a x +g(x) +\ln x -\ln(1/\theta_a) 
          \right\}\\ 
        &\ge 
        \frac{1}{4}
        \exp\{
          -\lambda_a^+ x +g(x)
          \}
          \ge \frac{1}{4}
          \exp \left\{
            \delta_1 g(x)
            \right\} 
           \ge  
           \frac{1}{4}
          \exp \left\{
            \delta_1 g(\varepsilon_a/a)
            \right\} 
            \to +\infty, 
          \end{align*}  
          \begin{comment}
      First we will show that 
      uniformly in $x\le (1-\delta)x(a)$, 
      \begin{equation}
        \label{exp-dominates-local}
      \frac{2}{\theta_a}
      \overline F(x\vee \varepsilon_a/a)
      =o(\theta_a e^{-\theta_a x}). 
      \end{equation}
      Equivalence \eqref{exp-dominates-local} clearly holds 
      uniformly in  $x\le 1/\theta_a$ 
      since $\overline F(\varepsilon_a/a) = o(a^2)$.
      Now for $x>1/\theta_a$ we have, 
      as $a\to 0$, 
      \begin{align*}
        \frac{\theta_a e^{-\theta_a x}}{\frac{1}{\theta_a}\overline F(x)}
        &\sim \frac{2a^2e^{-\theta_a x}}{x^{-2}e^{-g(x)}}
        =\exp
        \left\{
          -\theta_a x +g(x) +2\ln x -2\ln(1/\theta_a) 
          \right\}\\ 
        &\ge 
        \exp\{
          -\theta_a x +g(x)
          \}  
          =
          \exp\left\{
          -x\left(\theta_a  -\frac{g(x)}{x^{\gamma_0}} x^{\gamma_0-1}\right)
          \right\} \\
         &\ge 
         \exp\left\{
          -x\left(\theta_a  -\frac{g(x(a))}{x(a)^{\gamma_0}
          }((1-\delta)x(a))^{\gamma_0-1}\right)
          \right\}
          \\
          &=
          \exp\left\{
          -x\left(\theta_a  - 2a
          (1-\delta)^{\gamma_0-1}\right)
          \right\} \to \infty,  
          \end{align*}  
          uniformly in $x\le (1-\delta)x(a)$,  
          since $\theta_a\sim 2a$ and $(1-\delta)^{\gamma_0-1}>1$. 
        \end{comment}    
        uniformly in $x\le (1-\delta)x(a)$,  
          since $\varepsilon_a/a\to\infty$. 
         This proves  \eqref{exp-dominates}. 
      \end{proof}
        
    Next we will show that for $x\ge (1+\delta)x(a)$ 
    the subexponential term dominates 
    over the exponential term in $\overline G_\pm$. 
    \begin{lemma}
      \label{lem.subexp-domination}
      Let $\varepsilon_a\to 0$
      and let $\lambda^\pm_a\sim \theta_a\sim 2a/\sigma^2$ as $a\to 0$.
      Then, for any $\delta>0$, 
      uniformly in $x\ge (1+\delta)x(a)$, 
      \begin{equation}
        \label{eq.subexp-dominates}
        e^{-\lambda^\pm_a x}
        =
        o\left(
        \frac{2}{\theta_a }
        \overline F^I(x+ 2\varepsilon_a/a)
        \right). 
      \end{equation}  
    \end{lemma} 
    \begin{proof}
      To simplify notation we assume that $\sigma=1$. 
      In view of the insensitivity property \eqref{subexp-long-tails}
      it is sufficient to show that uniformly in $x>(1+\delta)x(a)$, 
      \begin{equation}
        \label{eq.subexp-dominates2}
        e^{-(\theta_a\pm c_a) x}
        =
        o\left(
        \frac{2}{\theta_a }
        \overline F^I(x)
        \right). 
      \end{equation}  
      Now 
      note that for $x>(1+\delta)x(a)$, 
      it follows from \eqref{eq.bounds.via.xa.below}
      \begin{align}
        \nonumber 
        -\lambda^\pm_a x +g(x) 
        &=- x \left(
        \lambda^\pm_a -\frac{g(x)}{x} 
        \right) \\ 
        &\le  
        \label{eq.minimizer.subex.0}
        -x
        \left(
        \lambda^\pm_a-\frac{g(x(a))}{x(a)} \left(\frac{x(a)}{x}\right)^{1-\gamma_0} 
        \right)\\ 
        &
        \label{eq.minimizer.subex}
        \le -x
        \left(
        \lambda^\pm_a -2a(1+\delta)^{1-\gamma_0}
        \right)
        \le -2\delta_1 \theta_a x,
  \end{align}  
      for some $\delta_1 \in (0,1-(1+\delta)^{1-\gamma_0}) $ 
      and all sufficiently small $a$, since 
      $-\lambda^\pm_a\sim \theta_a\sim 2a$.  
      Therefore, using the latter inequality  
      and the lower bound  \eqref{eq.bounds.integrated.tail.lower} 
      we obtain for a positive constant $c$, 
      uniformly in $x>(1+\delta)x(a)$, 
      as $a\to 0$, 
      \begin{align*}
        \frac{e^{-\lambda^\pm_a x}}{\frac{2}{\theta_a}\overline F^I(x)}
        &\le c 
        \frac{\theta_a xg(x) e^{-\lambda^\pm_a x}}{e^{-g(x)}}
        \le c \theta_a xg(x) e^{-2\delta_1 \theta_a x}\\ 
        &= c (\theta_a x)^2 \frac{g(x)}{\theta_a x }e^{-2\delta_1 \theta_a x}
        \le c (\theta_a x)^2 \frac{g(x(a))}{\theta_a x(a) }e^{-2\delta_1 \theta_a x}
        \le 3 c (\theta_a x)^2 e^{-2\delta_1 \theta_a x}\\ 
        &\le 3 (\theta_a (1+\delta)x(a))^2 
        e^{-2\delta_1 \theta_a (1+\delta)x(a)}
        \to 0, 
      \end{align*}  
      since $x^2e^{-x}$ is eventually decreasing to $0$ 
      and $ax(a)\to \infty$ by \eqref{eq.bound1.xa1}.
    \end{proof}      
    Next, we will analyse the convolution term 
    \begin{equation}
      \label{def.ixa}
      I(x,a):=
            \int_0^{x^+} 
            \left(
            \frac{1}{\lambda_a^\pm}
            +(x^+-z)
            \right)
            e^{-\lambda_a^\pm(x^+-z)}
            \overline F_+(z)
            dz.
    \end{equation}  
    \begin{lemma}
      \label{lem.simplify.ixa}
      Let $\varepsilon_a$ 
          and $c_a$ be defined according to Lemma~\ref{lem.sequences},  
          let $\theta_a\sim 2a/\sigma^2$, 
          as $a\to 0$. 
          Then, uniformly in $x$, as $a\to 0$,
          \begin{multline}
            \label{eq.simplifi.ixa}
            I(x,a) = (1+o(1)) 
            \int_{\varepsilon_a/a}^{x} 
        \left(
        \frac{1}{\theta_a}
        +(x-z)
        \right)
        e^{-\lambda_a^\pm(x-z)}
        \overline F(z)
        dz \mathbf{1}(x>\varepsilon_a/a)
        \\ +
        o\left(
          e^{-\theta_a x}+\frac{1}{\theta_a} \overline F^I(x)
        \right).
          \end{multline}  
    \end{lemma} 
    \begin{proof}
      By definition and change of variables in the integral, 
      \begin{align}
      I(x,a)
      \nonumber
      &= 
      \int_{0}^{x+\varepsilon_a/a} 
            \left(
            \frac{1}{\lambda_a^\pm}
            +x+\frac{\varepsilon_a}{a}-z
            \right)
            e^{-\lambda_a^\pm(x+\varepsilon_a/a-z)}
            \overline F(z+\varepsilon_a/a)dz\\
      &
      \label{lem19.eq0}
      =
            \int_{\varepsilon_a/a}^{x+2\varepsilon_a/a} 
            \left(
            \frac{1}{\lambda_a^\pm}
            +x+\frac{2\varepsilon_a}{a}-z
            \right)
            e^{-\lambda_a^\pm(x+2\varepsilon_a/a-z)}
            \overline F(z)dz. 
      \end{align}     
      The statement of Lemma clearly holds if 
      $x\le\varepsilon_a/a $ and therefore, we will consider 
      only $x>\varepsilon_a/a$. 
      Note that, as  
     $a\to 0$, 
     \begin{equation}
      \label{lem19.eq2}
      \int_x^{x+2\varepsilon_a/a}
      \left(
      \frac{1}{\lambda_a^\pm}
      +x+\frac{2\varepsilon_a}{a}-z
      \right)
      e^{-\lambda_a^\pm(x+2\varepsilon_a/a-z)}
            \overline F(z)dz
        \le \frac{2}{\theta_a} \frac{2\varepsilon_a}{a} \overline F(x).    
     \end{equation} 
     Then, for $x>\frac{1}{2}x(a)$, 
     using the lower bound 
     \eqref{eq.bounds.integrated.tail.lower}, 
     we obtain 
     \begin{align} 
      \nonumber
      \frac{\frac{2}{\theta_a} \frac{\varepsilon_a}{a} \overline F(x)}{\frac{2}{\theta_a} \overline F^I(x)}
      &\le C\frac{\frac{\varepsilon_a}{a} x^{-2}e^{-g(x)}}{x^{-1}g(x)^{-1}e^{-g(x)}}
      = C \frac{\varepsilon_a}{a} \frac{g(x)}{x}
      \le 
      C \frac{\varepsilon_a}{a} \frac{g(x(a)/2)}{x(a)/2}\\
      &\le 2C \frac{\varepsilon_a}{a} 
      \frac{g(x(a))}{x(a)}
      =\frac{2C \varepsilon_a}{\sigma^2} \to 0.
      \label{lem19.eq3}
     \end{align} 
     For $x\le \frac{1}{2}x(a)$, 
     we obtain, as $a\to 0$,  
     using \eqref{eq.minimizer}
     \begin{align} 
      \nonumber
      \frac{\frac{2}{\theta_a} \frac{\varepsilon_a}{a} \overline F(x)}{e^{-\theta_a x}}
      &\le 3\frac{\frac{\varepsilon_a}{a^2} x^{-2}e^{-g(x)}}{e^{-\theta_a x}}
      \le \frac{3}{\varepsilon_a} 
      e^{\theta_a x-g(x)}
      \le 
      \frac{3}{\varepsilon_a} 
      e^{-\delta_1 g(x)}
      \\
      \nonumber 
      &\le 
      \frac{3}{\varepsilon_a} 
      e^{-\delta_1 g(\varepsilon_a/a)}
      = 3 (ax(a))^{\delta_0}
      e^{-\delta_1 g(\varepsilon_a/a)}\\ 
      &\le 
      3
      \left(
      \frac{g(\varepsilon_a/a)}{2}
      \right)^{\frac{\delta_0}{1-\gamma_0(1+\delta_0)}}
      e^{-\delta_1 g(\varepsilon_a/a)}
      \to 0
      \label{lem19.eq4},   
     \end{align}
     where we used the definition \eqref{def.varepsilon.a} of 
     $\varepsilon_a$, the upper bound 
     \eqref{axa} and the fact that $\varepsilon_a/a\to +\infty$.

     Equations \eqref{lem19.eq2}, \eqref{lem19.eq3} 
     and \eqref{lem19.eq4} together 
     with \eqref{lem19.eq0} imply that 
     uniformly in $x$, 
     \begin{multline*}
      I(x,a) = 
      \int_{\varepsilon_a/a}^{x} 
            \left(
            \frac{1}{\lambda_a^\pm}
            +x+\frac{2\varepsilon_a}{a}-z
            \right)
            e^{-\lambda_a^\pm(x+2\varepsilon_a/a-z)}
            \overline F(z)dz \mathbf{1}(x>\varepsilon_a/a)
            \\ 
        +     o\left(
          e^{-\theta_a x}+\frac{1}{\theta_a} \overline F^I(x)
        \right). 
     \end{multline*}
     The statement of the Lemma now follows from the 
     observation that  
      $\lambda_a^\pm (2\varepsilon_a/a)\to 0$, as 
      $a\to 0$.  
    \end{proof}   

        \begin{lemma}
          \label{lem.conv-o}
          Let $\varepsilon_a$ 
          and $c_a$ be defined according to Lemma~\ref{lem.sequences},  
          let $\theta_a\sim 2a/\sigma^2$, 
          as $a\to 0$.
          Then, for any $\delta\in(0,1)$, 
          uniformly in $x\le (1-\delta)x(a)$, 
          \begin{equation}
            \label{conv-o}
            I(x,a)
            =o\left(
              e^{-\theta_a x}+\frac{1}{\theta_a} \overline F^I(x)
            \right), 
          \end{equation}  
          where $I(x,a)$ is defined in \eqref{def.ixa}.
        \end{lemma} 
        \begin{proof}
          Asymptotics equivalence clearly holds 
          for $x\le \varepsilon_a/a$ by Lemma~\ref{lem.simplify.ixa}. 
          Hence we will assume that $x>\varepsilon_a/a$.  
          Using \eqref{eq.simplifi.ixa} 
          and the insensitivity property 
           \eqref{exp-long-tails} we obtain, 
          uniformly in $x\le (1-\delta)x(a)$,  
          \begin{align*}
            I(x,a)
            &\sim \int_{\varepsilon_a/a}^{x} 
            \left(
            \frac{1}{\theta_a}
            +x-z
            \right)
            e^{-\theta_a(x-z)}
            \overline F(z)dz+
            o\left(
              e^{-\theta_a x}+\frac{1}{\theta_a} \overline F^I(x)
            \right)
            \\ 
            &\le
            2x(a)e^{-\theta_a x}
            \int_{\varepsilon_a/a}^{x} e^{\theta_a z} 
            \overline F(z) dz +o\left(
              e^{-\theta_a x}+\frac{1}{\theta_a} \overline F^I(x)
            \right)\\ 
            &\le 2x(a)\int_{\varepsilon_a/a}^{(1-\delta)x(a)} e^{\theta_a z} 
              \overline F(z)dz 
              +o\left(
                e^{-\theta_a x}+\frac{1}{\theta_a} \overline F^I(x)
              \right). 
        \end{align*}  
         Thus,  we are left to show that 
         \begin{equation}
          \label{eq.domin.last}
          x(a)\int_{\varepsilon_a/a}^{(1-\delta)x(a)} e^{\theta_a z} 
          \overline F(z)dz
        \to 0, \quad a\to 0. 
         \end{equation}
         %It immediately follows from 
         %$\overline F(\varepsilon_a/a) = o(a/x(a))$
         %that 
         %\begin{equation}
         % \label{nakonets2}
         %x(a)\int_{\varepsilon_a/a}^{1/a} e^{\theta_a z} 
         % \overline F(z)dz
         % \le e^{\theta_a/a} \overline F(\varepsilon_a/a) 
         % \frac{x(a)}{a}\to 0.
         %\end{equation}
         %Next, 
         It follows from \eqref{eq.minimizer} 
         that for some $\delta_1>0$, 
         \begin{align}
          \nonumber
          x(a)\int_{\varepsilon_a/a}^{(1-\delta)x(a)} 
          e^{\theta_a z} \overline F(z)dz
          &\sim x(a) 
          \int_{\varepsilon_a/a}^{(1-\delta)x(a)}
          e^{\theta_a z - g(z)}z^{-2}dz \\ 
          \nonumber
          &\le x(a) 
          \int_{\varepsilon_a/a}^{(1-\delta)x(a)}
          e^{-\delta_1 g(z)}z^{-2}dz\\
          &
          \nonumber
          \le 
          x(a) e^{-\delta_1g(\varepsilon_a/a)}
          \int_{\varepsilon_a/a}^\infty z^{-2}dz \\ 
          \label{domin.second.integral}
          &= \frac{a x(a)}{\varepsilon_a} e^{-\delta_1g(\varepsilon_a/a)}.
          %&\le 
          %\frac{ax(a)}{\varepsilon_a}
          %e^{-\delta_1 g(\varepsilon_a/a)}.  
         \end{align} 
         Using the definition \eqref{def.varepsilon.a} of $\varepsilon_a$ 
         and \eqref{axa} we hence obtain 
           \begin{equation}
            \label{nakonets}
            x(a)\int_{\varepsilon_a/a}^{(1-\delta)x(a)} 
            e^{\theta_a z} \overline F(z)dz
            \le 
            \left(\frac{g(\varepsilon_a/a)}{2}\right)^{\frac{1+\delta_0}{1-\gamma_0(1+\delta_0)}}
            e^{-\delta_0g(\varepsilon_a/a)}\to 0, a\to 0. 
     \end{equation} 
     Equation \eqref{nakonets} 
     implies \eqref{eq.domin.last} and hence 
     \eqref{conv-o}.
    \end{proof}

    \begin{lemma}
      \label{lem.conv.2}
      Let $\varepsilon_a$ 
      and $c_a$ be defined according to Lemma~\ref{lem.sequences},  
      let $\theta_a\sim 2a/\sigma^2$, 
      as $a\to 0$.
      Let $A_a\uparrow \infty$, as $a\to 0$,  be an arbitrary function.  
      Then, for any $\delta\in(0,1)$, 
      uniformly in $x\in [(1-\delta)x(a), A_ax(a)]$, 
      \begin{equation}
        \label{conv-o-2}
        I(x,a)\sim 
        \int_{(1-\delta)x(a)}^{x} 
        \left(
        \frac{1}{\theta_a}
        +(x-z)
        \right)
        e^{-\theta_a(x-z)}
        \overline F(z)
        dz
        +
        o\left(
          e^{-\theta_a x}+\frac{1}{\theta_a} \overline F^I(x)
        \right), 
      \end{equation}  
      and uniformly in $x\ge A_ax(a)$, 
      \begin{equation}
        \label{conv-o-3}
        I(x,a)=
        o\left(
          e^{-\theta_a x}+\frac{1}{\theta_a} \overline F^I(x)
        \right), 
      \end{equation}
      where $I(x,a)$ was defined in \eqref{conv-o}. 
    \end{lemma}  
    \begin{proof}
     First we split the integral in \eqref{eq.simplifi.ixa} 
     to obtain  
     \begin{align*}
      I(x,a) &\sim  
      \left(
        \int_{\varepsilon_a/a}^{(1-\delta)x(a)} 
        +
        \int_{(1-\delta)x(a)}^x 
      \right)
            \left(
            \frac{1}{\theta_a}
            +x-z
            \right)
            e^{-\lambda_a^\pm(x-z)}
            \overline F(z)dz 
            + o\left(
              e^{-\theta_a x}+\frac{1}{\theta_a} \overline F^I(x)
            \right)\\
            &:=I_1(x,a)+I_2(x,a)+o\left(
              e^{-\theta_a x}+\frac{1}{\theta_a} \overline F^I(x)
            \right). 
     \end{align*} 
     Similarly to Lemma~\ref{lem.conv-o}, using 
     \eqref{eq.minimizer}, 
     \begin{align*}
      I_1(x,a)&\le 2x e^{-\lambda_a^+ x} 
      \int_{\varepsilon_a/a}^{(1-\delta)x(a)}e^{\lambda^+_a z}\overline F(z)dz\\ 
      &\sim 2x e^{-\lambda^+x} 
      \int_{\varepsilon_a/a}^{(1-\delta)x(a)} 
      e^{\lambda^+_a z-g(z)}z^{-2}dz 
      \le 2x e^{-\lambda^+x} 
      \int_{\varepsilon_a/a}^{(1-\delta)x(a)} 
      e^{-\delta_1 z}z^{-2}dz\\ 
      &\le e^{-\lambda_a^+ x}  \frac{2a x} {\varepsilon_a} e^{-\delta_1 g(\varepsilon_a/a)}.
     \end{align*} 
     Using \eqref{nakonets} we obtain, uniformly in $x\ge (1-\delta)x(a)$
     \begin{equation}
      I_1(x,a)\le 
      (1+o(1)) e^{-\lambda_a^+ x} \frac{2x}{x(a)} 
      \left(\frac{g(\varepsilon_a/a)}{2}\right)^{\frac{1+\delta_0}{1-\gamma_0(1+\delta_0)}}
            e^{-\delta_1g(\varepsilon_a/a)}.
     \end{equation} 
     Then, there exists $\widetilde A_a\uparrow \infty$ such that 
     %$$
     %A_a = o\left(
     % \left(\frac{g(\varepsilon_a/a)}{2}\right)^{\frac{1+\delta_0}{1-\gamma_0(1+\delta_0)}}
     % e^{-\delta_0g(\varepsilon_a/a)}
     %\right)
     %$$
     we obtain uniformly in $x\in [(1-\delta)x(a), \widetilde A_a x(a)]$ that 
    \begin{equation}
      \label{i1xa}
    I_1(x,a) = o(e^{-\theta_a x}).  
    \end{equation}
    Next, using the lower bound \eqref{eq.bounds.integrated.tail.lower} 
    and then \eqref{eq.minimizer.subex},  
    we obtain uniformly in $x \ge \widetilde A_a x(a)$, 
    \begin{align*}
      \frac{I_1(x,a)}{\frac{2}{\theta_a}\overline F^I(x)}
      &\le o(1) \frac
      {x e^{-\lambda^+ x}}
      {\frac{2}{\theta_a} x^{-1}g(x)^{-1}e^{-g(x)}} 
      \le o(1) x^2 g(x) 
      e^{-2\delta_1 \theta_a x}.
    \end{align*}  
    Now note that by \eqref{eq.bounds.via.xa.below} 
    for $x\ge \widetilde A_a x(a)$ 
    \begin{align*}
      \theta_a x  
      \ge 
      g(x)\frac{\theta_a}{2a}
      \widetilde A_a^{1-\gamma_0}
      \ge \frac{1}{2}g(x) \widetilde A_a^{1-\gamma_0}. 
    \end{align*}  
    Therefore, uniformly in $x \ge \widetilde A_a x(a)$, 
    \begin{align*}
      \frac{I_1(x,a)}{\frac{2}{\theta_a}\overline F^I(x)}
      &\le o(1) x^{2} g(x) e^{-\widetilde A_a^{1-\gamma_0} \delta_a g(x)} 
      \to 0,
    \end{align*}  
    using the facts that $g(x)\ge \varepsilon \ln x$ for some $\varepsilon>0$ 
    and that $\widetilde A_a\to \infty$. The proves \eqref{conv-o-2}. 

    Now consider the case  $x\ge A_ax(a)$. 
    Here, we are left to show that uniformly in $x\ge A_ax(a)$, 
    \begin{equation}
      \label{i2xa}
      I_2(x,a) = o\left(
        \frac{2}{\theta_a} \overline F^I(x)
        \right). 
    \end{equation}
    %Let $C>1$ be a fixed constant and split 
    Split 
    \begin{align*}
      I_2(x,a) &= 
      2\left(
        \int_{(1-\delta)x(a)}^{x/2}
        +\int_{x/2}^x 
      \right)
      \left(
        \frac{1}{\theta_a} + x-z 
      \right)
      e^{-\lambda^\pm (x-z)} 
      \overline F(z)dz \\
      &:=I_{21}(x,a)+I_{22}(x,a). 
    \end{align*}  
    \begin{comment}
    We will use the following bound 
    for $x\ge C x(a)$, 
    \begin{align*}
      -\lambda^\pm_a x+g(x) 
      &\le g(x) \left(
      -\theta_a \frac{x}{g(x)}+1 
      \right)
      \le 
      g(x) \left(
      -\theta_a \frac{Cx(a)}{g(Cx(a))}+1 
      \right)\\ 
      &\le 
      g(x) \left(
      -\theta_a (Cx(a))^{1-\gamma_0}\frac{x(a)^{\gamma_0}}{g(x(a))}+1 
      \right)\\ 
      &\le g(x) 
      \left(
        -\frac{\theta_a}{2a} C^{1-\gamma_0}+1
      \right)
      \le -\frac{C^{1-\gamma_0}}{2}g(x).
    \end{align*}  
  \end{comment}  
  First, 
    \begin{align*}
      I_{21}(x,a)
      \le 2x e^{-\lambda^\pm_a x/2}
      \int_{x(a)/2}^\infty \overline F(z) dz 
      =2xe^{-\lambda^\pm_a x/2} o(1/x(a)).
    \end{align*}  
    Therefore, using \eqref{eq.minimizer.subex.0}, 
    as $a\to 0$, 
    \begin{align*}
      \frac{I_{21}(x,a)}{
        \frac{1}{\theta_a} 
        \overline F^I(x)
      }
      &\le o(\theta_a /x(a))
      x^2 g(x) e^{-\lambda^\pm x/2 +g(x)}
      \le 
      o(\theta_a /x(a))
      x^2 g(x)e^{-x(\lambda^\pm_a/2 +2a A_a^{\gamma_0-1})}\\ 
      &\le o(1) (\theta_a x)^2 g(x) e^{-\theta_a x/4}
      =o(1) (\theta_a x)^2 x \frac{g(x)}{x} e^{-\theta_a x/4}\\ 
      &\le o(1) (\theta_a x)^3 e^{-\theta_a x/4} 
      \le o(1) (\theta_a A_a x(a))^3 e^{-\theta_a A_a x(a)/4} 
      \to 0,
    \end{align*}  
    since $\theta_a x(a)\to 0$.
  Next, 
    \begin{align*}
      I_2(x,a)
      &= 
      2\int_{0}^{x/2} 
      \left(
        \frac{1}{\theta_a} + z 
      \right)
      e^{-\lambda^\pm z} 
      \overline F(x-z)dz  \\
      &\sim 
      2e^{-g(x)}\int_{0}^{x/2} 
      \left(
        \frac{1}{\theta_a} + z 
      \right)
      (x-z)^{-2}
      e^{-\lambda^\pm z-g(x-z)+g(x)} 
      dz\\ 
      &\le 
      8e^{-g(x)}x^{-2}
      \int_{0}^{x/2} 
      \left(
        \frac{1}{\theta_a} + z 
      \right)
      e^{
        -\lambda^\pm z 
        + \gamma_0 z \frac{g(A_a x(a)/2)}{A_ax(a)/2}
      }dz \\ 
      &\le 
      8e^{-g(x)}x^{-2}
      \int_{0}^{x/2} 
      \left(
        \frac{1}{\theta_a} + z 
      \right)
      e^{-\theta_a z/2}
      dz
      \le 
      \frac{50}{\theta_a^2} e^{-g(x)}x^{-2}.
   \end{align*}    
    Then, uniformly in $x\ge A_a x(a)$, 
    \begin{align*}
      \frac{I_2(x,a)}{\frac{2}{\theta_a}\overline F^I(x)}
      &\le 
      \frac{I_2(x,a)}{\frac{2}{\theta_a}(xg(x))^{-1}e^{-g(x)}}
      \le 
      \frac{\frac{50}{\theta_a^2} e^{-g(x)}{x^{-2}}}{\frac{2}{\theta_a}(xg(x))^{-1}e^{-g(x)}}\\ 
      &\le \frac{50}{\theta_a} \frac{g(x)}{x} 
      \le \frac{50}{\theta_a} \frac{g(A_ax(a))}{Ax(a)}
      \le  \frac{50}{\theta_a} 
      A_a^{\gamma_0-1}
      \frac{g(x(a))}{x(a)}
      \le 
      \frac{100 a}{\theta_a} A_a^{\gamma_0-1} \to 0. 
    \end{align*}  
    This completes the proof. 
      \end{proof}  

      We are now in position to prove Theorem~\ref{thm.main}. 
      \begin{proof} (of Theorem~\ref{thm.main})
        Let $x(a)$ be defined according to \eqref{eq.boundary}. 
        Let $c_a$ and $\varepsilon_a$ 
        be defined as in Lemma~\ref{lem.sequences}. 
        For this choice of $c_a$ and , $\theta_a$ solving \eqref{def.theta.xa} 
        will be solving \eqref{def.theta} as well. 
        Also, by Lemma~\ref{lem.sequences}, 
         $F$ will satisfy \eqref{subex} and $\overline F(\varepsilon_a/a)$. 
        Hence, all conditions of Theorem~\ref{prop1} are met and 
        for arbitrary $\alpha$ there exist $a_0$ such 
        that for $a\in(0,a_0)$ 
        the lower  bound \eqref{eq.lower.bound} and 
        the upper bounds \eqref{eq.upper.bound} hold. 
        Since $\alpha>0$ is arbitrary, 
        we can put $\alpha=0$ in the definition of  
         $\overline G_\pm$ and  
        the lower and upper bounds will still 
        hold asymptotically as $a\to 0$, uniformly in $x$.

        By Lemma~\ref{lem.exp-domination} 
        and Lemma~\ref{lem.conv-o}, 
        uniformly in  $x\le (1-\delta)x(a)$, 
        $$
        \overline G_+(x)\sim \overline G_-(x)
        \sim e^{-\theta_a x},
        $$
        as $a\to 0$.  This 
        implies 
        that  uniformly in  $x\le (1-\delta)x(a)$, 
        $$
        \pr(M^{(a)}>x) \sim e^{-\theta_a x},\quad  a \to 0. 
        $$

        Next consider $x\ge A_ax(a)$.  In this case, 
        applying Lemma~\ref{lem.subexp-domination}, 
        insensitivity Lemma~\ref{lem.insensitivity} 
        and \eqref{conv-o-3} we obtain 
        that uniformly in 
        $x\ge A_ax(a)$
        $$
        \overline G_+(x)\sim \overline G_-(x)
        \sim 
        \frac{2}{\theta_a}
        \overline F^I(x)
        \sim 
        \frac{1}{a}
        \overline F^I(x)
        ,
        $$
        as $a\to 0.$ 
        This implies 
        that  uniformly in $x\ge A_ax(a)$, 
        $$
        \pr(M^{(a)}>x) \sim \frac{1}{a}\overline F^I(x),\quad  a \to 0. 
        $$

        Finally, for 
        $x\in [(1-\delta)x(a), A_ax(a)]$  the result follows 
        from \eqref{conv-o-2} and the insensitivity 
        Lemma~\ref{lem.insensitivity}.

      \end{proof}

      \section{Solution to equation \eqref{def.theta.xa}}
    \label{sec.sol.theta.xa}
    \begin{lemma}
      \label{equation:rv}
      Let $\theta_a=2a/\sigma^2$. 
      Assume that $\e|X|^\gamma<\infty$ 
       for some  $\gamma\in [2,3)$. 
      Then,  
      $$
      \e[e^{\theta_a X^{(a)}}; X^{(a)}\le 1/a ] 
      =1+O(a^\gamma). 
      $$
    \end{lemma}  
    \begin{proof}
      We have, 
      $$
        \e[e^{\theta_a X^{(a)}}; X^{(a)}\le 1/a]=
        \e[e^{\theta_a X^{(a)}}; |X^{(a)}|\le 1/a] 
        +\e[e^{\theta_a X^{(a)}}; X^{(a)}\le -1/a]\\  
      $$
      Put $R(x) = e^x-1-x-x^2/2$.  
      Then, it follows from the Taylor formula that 
      $|R(x)|\le 3|x|^3$ for $x: |x|\le 1$. 
      Hence, 
      \begin{align*}
      &\e[e^{\theta_a X^{(a)}}; |X^{(a)}|\le 1/a] 
      =\e\left[1+\theta_a X^{(a)} 
      +\frac{1}{2}(\theta_a X^{(a)})^2;
      |X^{(a)}|\le 1/a] 
       \right]\\ 
      &\hspace{1cm}+\e[R(\theta_a X^{(a)});|X^{(a)}|\le 1/a]\\ 
      &=  
      \e\left[1+\theta_a X^{(a)} 
      +\frac{1}{2}(\theta_a X^{(a)})^2
      \right] 
      -\pr(|X^{(a)}|>1/a)
      -\theta_a\e[X^{(a)};|X^{(a)}|>1/a]\\
      &\hspace{1cm}-
      \frac{\theta_a^2}{2}\e[(X^{(a)})^2;|X^{(a)}|>1/a]
      +\e[R(\theta_a X^{(a)});|X^{(a)}|\le 1/a].
      \end{align*}
    Now, by the Markov inequality, 
    \begin{align*}
    &\pr(|X^{(a)}|>1/a) +
    \theta_a\e[X^{(a)};|X^{(a)}|>1/a] 
    +\frac{\theta_a^2}{2}\e[(X^{(a)})^2;|X^{(a)}|>1/a]\\
    &\hspace{1cm}+\frac{\theta_a^3}{2}\e[|X^{(a)}|^3;|X^{(a)}|>1/a]
    \le C a^\gamma \e[(X^{(a)})^\gamma] = O(a^\gamma)
    \end{align*}
    and 
    $$
    \e[e^{\theta_a X^{(a)}}: X^{(a)}\le -1/a]
    \le \pr(X^{(a)}\le -1/a) 
    \le a^\gamma \e[(X^{(a)})^\gamma] = O(a^\gamma).
    $$
    Hence, 
    \begin{align*}
    \e[e^{\theta_a X^{(a)}}; X^{(a)}\le 1/a ] 
    &=
    \e\left[1+\theta_a X^{(a)} 
      +\frac{1}{2}(\theta_a X^{(a)})^2
      \right] +O(a^\gamma)\\
    &=1-a\theta_a+\frac{\theta_a^2}{2}(\sigma^2+a^2)
    +O(a^\gamma) 
    =1+O(a^\gamma), 
    \end{align*}
    as required. 
  \end{proof}  
  Accuracy of the solution to the Cram\'er equation \eqref{def.theta.xa} 
  given by Lemma~\ref{equation:rv} is sufficient for $\gamma_0<1/2$. 
  For $\gamma_0\in[1/2,1)$ we can construct recursively polynomial 
  approximation which will give sufficient accuracy. 
  This construction will be described in the following lemma
  \begin{lemma}
    \label{lem.eq.weibull}
    Assume that $\e|X|^\gamma<\infty$ for some $\gamma\in [n+1,n+2).$
    Then there exists a polynmial of degree $n$  
    $$
    \theta_a^{(n)} = \sum_{k=1}^n C_k a^n 
    $$
    such that 
    $$
    \e[e^{\theta_a^{(n)}X^{(a)}};X^{(a)}\le 1/a] = 1+O(a^\gamma), 
    \quad a\to 0.
    $$
  \end{lemma}  
  \begin{proof}
    Put 
    $$
    \mu_k^{(a)} = \e[(X^{(a)})^k] =  \e[(X-a)^k].
    $$
    Clearly $\mu_k^{(a)}$ is a polynomial of degree $k$ whose 
    coefficients are defined by first $k$ moments. 
    We will first construct inductively polynomials such that 
    \begin{equation}
      \label{lem23.eq1}
      \sum_{k=1}^{n+1} \frac{1}{k!} (\theta_a^{(n)})^k \mu_k^{(a)} = O(a^{n+2}).
    \end{equation}    
    For $n=1$ we put $C_1=2/\sigma^2$ and thus 
    $\theta_a^{(1)}=\frac{2}{\sigma^2}a$. 
    Clearly for  $\theta_a^{(1)}$ equation \eqref{lem23.eq1} holds 
    since 
    $$
    \sum_{k=1}^{2} \frac{1}{k!} (\theta_a^{(1)})^k \mu_k^{(a)}
    =\frac{2a}{\sigma^2}(-a) +\frac{1}{2}
    \left(\frac{2a}{\sigma^2}\right)^2(\sigma^2+a^2)
    =\frac{2a^4}{\sigma^4}
    =O(a^3).
    $$
    Now suppose that we have constructed $\theta_a^{(n)}$ 
    and we will construct $\theta_a^{(n+1)}$ 
    satisfying \eqref{lem23.eq1} for $n+1$. 
    For that we  put  
    $$
    \theta_a^{(n+1)}  = \theta_a^{(n)}+C_{n+1} a^{n+1}
    $$
    and will be looking for a suitable $C_{n+1}$. 
     Since $\theta_a^{(n)} \sim 2a/\sigma^2$, 
    \begin{align*}
      &\sum_{k=1}^{n+2} \frac{1}{k!} (\theta_a^{(n)}+C_{n+1}a^{n+1})^k \mu_k^{(a)}
      =
      \sum_{k=1}^{n+2} \frac{1}{k!} (\theta_a^{(n)})^k \mu_k^{(a)}\\ 
      &\hspace{1cm} +C_{n+1}a^{n+1}(-a)
      +C_{n+1} \theta_a^{(n)} a^{n+1}\sigma^2
      +O(a^{n+3}) \\ 
      =
      &\sum_{k=1}^{n+2} \frac{1}{k!} (\theta_a^{(n)})^k \mu_k^{(a)}
      +C_{n+1} a^{n+2} +O(a^{n+3}) .  
    \end{align*}
    Since  both  $\theta_a^{(n)}$ 
    and $\mu_k^{(a)}$ are polynomials in $a$ 
    for the induction assumption \eqref{lem23.eq1} to hold, 
    $$
    \sum_{k=1}^{n+2} \frac{1}{k!} (\theta_a^{(n)})^k \mu_k^{(a)}
    =B_{n+2} a^{n+2} + O(a^{n+3}). 
    $$
    Therefore, we can simply put $C_{n+1}=-B_{n+2}$ 
    to obtain \eqref{lem23.eq1} for $n+1$. 
    It is clear from this construction that 
    $C_{n}$ depends only on first $n+1$ moments of $X$. 

    Once we have constructed the polynomial $\theta_a^{(n)}$ 
    we can proceed to the proof of the statement. 
    Put 
    $$
    R_{n+1}(x) = e^{x}-\sum_{j=0}^{n+1}\frac{x^j}{j!} 
    $$
    By Taylor's formula for $x$ such that $|x|\le  1$, 
    $|R_{n+1}(x)| \le \alpha_{n+1}|x|^{n+2} $ for 
    some constants $\alpha_{n+1}$. 
    First,
    \begin{align*}
      \e[e^{\theta_a^{(n)} X^{(a)}}; X^{(a)}\le 1/a]
      &=
      \e[e^{\theta_a^{(n)} X^{(a)}}; |X^{(a)}|\le 1/a]
      +
      \e[e^{\theta_a^{(n)} X^{(a)}}; X^{(a)}<-1/a]\\ 
      &=E_1(a)+E_2(a)
    \end{align*}  
    By the Markov inequality we immediately obtain, 
    as $a\to 0,$
    \begin{equation}
      \label{lem23.e2}
    E_2(a)
    \le \pr(X^{(a)}<-1/a) 
    \le  a^{\gamma}\e[(X^{(a)})^{\gamma};X^{(a)}<-1/a] 
    =o(a^\gamma),
    \end{equation}
    For the first summand we will apply the Taylor expansion, 
    \begin{align*}
      E_1(a)&=
      \e\left[
        \sum_{j=0}^{n+1} 
        \frac{(\theta_a^{(n)} X^{(a)})^j}{j!}
        ; |X^{(a)}| \le 1/a
      \right]
      +
      \e\left[
        R_{n+1}(\theta_a^{(n)} X^{(a)}) ; |X^{(a)}| \le 1/a
        \right]\\
     &= 
     \e\left[
        \sum_{j=0}^{n+1} 
        \frac{(\theta_a^{(n)} X^{(a)})^j}{j!}
      \right]
      -
     \e\left[
        \sum_{j=0}^{n+1} 
        \frac{(\theta_a^{(n)} X^{(a)})^j}{j!}
        ; |X^{(a)}| > 1/a
      \right]  \\ 
      &\hspace{1cm}+
      \e\left[
        R_{n+1}(\theta_a^{(n)} X^{(a)}) ; |X^{(a)}| \le 1/a
        \right]
        :=E_{11}(a)+E_{12}(a)+E_{13}(a).
    \end{align*}  
    By the defining property \eqref{lem23.eq1} 
    of $\theta_a^{(n)}$ we immediately obtain 
    \begin{equation}
      \label{lem23.e11}
    E_{11}(a) =  1+ O(a^{n+2}),\quad a\to 0.
    \end{equation}
Since $\theta_a^{(n)}\sim 2a/\sigma^2$ 
  and the family of random variables 
  $\{(X^{(a)})^{\gamma}\}_{a>0}$ is uniformly integrable, 
  \begin{align}
    \nonumber
    |E_{12}(a)|&\le 
    C \sum_{j=0}^{n+1}
      a^j \e[(X^{(a)})^j;|X^{(a)}|>1/a]\\ 
      &\le 
      C \sum_{j=0}^{n+1}
      a^\gamma \e[(X^{(a)})^\gamma;|X^{(a)}|>1/a]
      =o(a^\gamma).
      \label{lem23.e12}
  \end{align}  
  Finally, 
  \begin{align}
    \nonumber
    |E_{13}(a)|&\le \alpha_{n+1} 
    \e[|\theta_a^{(n)} X^{(a)}|^{n+2};|X^{(a)}|\le 1/a]
    \le C a^{n+2}
    \e[|X^{(a)}|^{n+2};|X^{(a)}|\le 1/a]
    \\ 
    &\le C a^{\gamma} \e[|X^{(a)}|^{\gamma}]
    =O(a^\gamma).
    \label{lem23.e13}
  \end{align}  
  Now the statement follows from \eqref{lem23.e2} -- \eqref{lem23.e13}. 
\end{proof} 
  
    \section{Proofs of Corollary~\ref{thm.rv} and  Corollary~\ref{thm.weibull}}
    \label{sec.final}
      To analyse the case $g'(x) = o(g(x)/x)$ we first   
      slightly improve the lower bound. 
    \begin{lemma}
        \label{lem.bounds.integrated.tail2}
        Let $F$ satisfy \eqref{sc1} and \eqref{sc2}. 
        Assume, in addition, that $g'(x) =o(g(x)/x)$. 
        Then, there exists a function $C(x) \uparrow \infty$
        such that 
        for any $\delta>0$
        there exists $x_0$  such that 
        for $x>x_0,$
        \begin{align}
          \label{eq.bounds.integrated.tail.lower2}
          \overline F^I(x)&
          %\ge 
          %c \gamma_0  (x^{2}g'(x))^{-1} e^{-g(x)}
          \ge C(x) (xg(x))^{-1}e^{-g(x)}.
        \end{align}  
      \end{lemma}

      \begin{proof}
        Let $C(x)\uparrow \infty$ be an increasing non-negative function such that
        \begin{equation}
          \label{def_cx}
          \sup_{y\ge x}\frac{yg'(y)}{g(y)}\ge \frac{1}{C(x)}.
        \end{equation}  
        Note that, as $x\to\infty$,  
       \begin{align}
        \nonumber
        \overline F^I(x) &\ge 
        \int_x^{x+C(x)x/g(x)}\overline F(y) dy 
        \sim\int_x^{x+C(x) x/g(x)}y^{-2}e^{-g(y)}dy \\ 
        \label{eq.lower.bound.11}
        &\ge e^{-g(x+C(x) x/g(x))} \int_x^{x+C(x) x/g(x)}y^{-2}dy
        \sim e^{-g(x+C(x)x/g(x))} \frac{C(x)}{xg(x)}
        \end{align} 
        Also, 
        \begin{align*}
          g(x+C(x)x/g(x)) - g(x)
          &=\int_{x}^{x+C(x) x/g(x)} 
          \frac{zg'(z)}{g(z)} \frac{g(z)}{z}dz \\ 
          &\le 
          \frac{1}{C(x)}\int_{x}^{x+xC(x)/g(x)} 
          \frac{g(z)}{z}dz 
          \le 1.
        \end{align*}  
        Plugging in the latter inequality in  \eqref{eq.lower.bound.11}
        we obtain  
        $$
        \overline F^I(x) 
        \ge (1+o(1)) 
        \frac{C(x)}{xg(x)} e^{-g(x)} e^{-1},
        $$
        which implies  \eqref{eq.bounds.integrated.tail.lower2}.
      \end{proof}    

  \begin{proof} (of Corollary~\ref{thm.rv})
    Conditions of the theorem imply that the 
    statement of Theorem~\ref{thm.main} holds. 
    It remains to show that  (i)  
    there exists $A_a\uparrow \infty$ such that 
    the convolution term in \eqref{eq.main} is negligible, that is 
    $$
      I_2(x,a):=2\int_{(1-\delta)x(a)}^x 
    \left(
    \frac{1}{\theta_a}
    +(x-z)
    \right) e^{-\theta_a(x-z)} \overline F(z)dz
    =\left(
      e^{-\theta_a x} +\frac{1}{\theta_a} \overline F^I(x)
    \right),
    $$
    as $a\to 0$,  
    and  (ii) that  $\theta=\frac{2a}{\sigma^2}$ satisfies \eqref{def.theta.xa}. 

    First we will show that 
    there exists $A_a\uparrow \infty$ such that 
    the convolution term in \eqref{eq.main} is negligible
    Let $C(x)$ be the function defined in \eqref{lem.bounds.integrated.tail2}
    and let 
    $A_a\uparrow\infty$ be such that 
    \begin{equation}
      \label{eq.aa}
    A_a^2 = o\left(
      \frac{1}{C((1-\delta)x(a))}
      \right), \quad a\to 0.
    \end{equation}
  We have, as $a\to 0$, 
    \begin{align*}
      I_2(x,a)
    &\le 2 
    \int_{(1-\delta)x(a)}^x 
    (x-z)
    e^{-\theta_a(x-z)} \overline F(z)dz
    \sim 
    2 
    \int_{(1-\delta)x(a)}^x 
    (x-z)
    e^{-\theta_a(x-z)} z^{-2}e^{-g(z)}dz\\ 
     &\le 
     \frac{2}
     {(1-\delta)^2} 
     \frac{1}{x(a)^2}
     e^{-g(x)}
     \int_{(1-\delta)x(a)}^x 
     (x-z)
    e^{-\theta_a(x-z)} e^{g(x)-g(z)}dz  
    \end{align*} 
    Now note that since $g'(x) = o(g(x/x))$ 
    there exists $\delta(a)\to 0$ such that 
    for $x\ge z\ge (1-\delta)x(a)$
    \begin{align*}
      g(x)-g(z) &= \int_z^x g'(t) dt 
      \le \delta(a)
      \int_z^x \frac{g(t)}{t}
      \le \delta(a)(x-z) \frac{g(1-\delta)x(a)}{(1-\delta)x(a)}\\ 
      &\le \frac{\delta(a) }{(1-\delta)} 
      \frac{g(x(a))}{x(a)} (x-z)
      = \frac{\delta(a) 2a}{(1-\delta)} 
        (x-z).
    \end{align*}      
    Since, $\delta(a) \to 0$ and $\theta_a\sim 2a$ 
    we obtain the following estimate, as $a\to 0$, 
    \begin{align*}
      I_2(x,a)
      &\le 
      \frac{2}
     {(1-\delta)^2} 
     \frac{1}{x(a)^2}
     e^{-g(x)}
     \int_{(1-\delta)x(a)}^x 
     (x-z)
    e^{-\theta_a(x-z)/2}dz 
    \le 
    \frac{8}
     {(1-\delta)^2\theta_a^2} 
     \frac{1}{x(a)^2}
     e^{-g(x)}
    \end{align*}
    Using the lower bound \eqref{eq.bounds.integrated.tail.lower2} 
    we obtain, 
    \begin{align*}
      \frac{I_2(x,a)}{\frac{2}{\theta_a}\overline F^I(x)}
      &\le 
      \frac{I_2(x,a)}{\frac{2}{\theta_a}C(x) 
      (xg(x))^{-1}e^{-g(x)}}
      \le 
      \frac{4}{\theta_a (1-\delta)^2}
      \frac{xg(x)}{x(a)^2 C((1-\delta)x(a))}\\ 
      &\le
      \frac{2}{(1-\delta)^2}
      \left(\frac{x}{x(a)}\right)^2
      \frac{g(x)}{x \theta_a} \frac{1}{C(1-\delta)x(a)}
      \le
      \frac{2}{(1-\delta)^2}
      \frac{2a}{\theta_a} \frac{A_a^2}{C(1-\delta)x(a)}
      \to 0, 
    \end{align*}    
    using \eqref{eq.aa}.  This proves that the convolution 
    term is negligible.

    Second, let $\theta=2a/\sigma^2$ and this choice will satisfy \eqref{def.theta.xa}. 
    For that note that $g'(x) =o(g(x)/x)$ implies that 
    $g(x)=o(x^{\delta_2}), x\to \infty$ for any $\delta_2>0$. 
    Then, \eqref{eq.boundary} implies that 
    $x(a) = o(a^{-\delta_2-1})$ for any $\delta_2>0$, as  $a\to 0$.
    Next,  by our assumptions  $\e[|X|^{2+\varepsilon}]<\infty$ 
    for some $\varepsilon\in(0,1)$. Then, by Lemma~\ref{equation:rv} 
    \begin{align*}
      \e[e^{\frac{2a}{\sigma^2}X^{(a)}};X^{(a)}\le 1/a] 
      =1+O(a^{2+\varepsilon}) = 
      1+o(a/x(a)),
    \end{align*}  
    as $x(a) = o(a^{-1-\varepsilon})$. 
  \end{proof}  
  
  \begin{proof} (of Corollary~\ref{thm.weibull})
    Conditions of the theorem imply that the 
    statement of Theorem~\ref{thm.main} holds. 
    Hence, it is sufficient to find the asymptotics for 
    the convolution term in \eqref{eq.main} for some 
    function $A_a\uparrow \infty$.

    It follows from the Karamata theorem that 
    $$
    g'(x) \sim \beta \frac{g(x)}{x},\quad x\to \infty. 
    $$
    Then, $\overline F$ satisfies \eqref{sc1} and \eqref{sc2} 
    for some  $\gamma_0\in(\beta, 1)$.  It also follows from the L'Hopital rule that 
    $$
    \overline F^I(x) \sim \frac{\overline F(x)}{g'(x)}, \quad x\to\infty. 
    $$
    Let $\delta$ be such that $(1-\delta)>\beta$.  
    Clearly, all conditions of Theorem~\ref{thm.main} 
    are met and we are left to find the 
    asymptotics of the convolution term. 
    
    If $x\in [(1-\delta)x(a), (1-\delta/2)x(a)$
    then by Lemma~\ref{lem.conv-o} the convolution term is negligible.
    It is not difficult to see  that the the third term in \eqref{eq.main.weibull} is of smaller 
    order than th e sum of the first and the second. Hence, the statement is 
    true for $x\in [(1-\delta)x(a), (1-\delta/2)x(a)$ and we will 
    consider only $x>(1-\delta/2)x(a)$.  

    Next,  for $x\ge (1-\delta)x(a)$,  any fixed constant $C$ 
    uniformly in $z$ such that $z \le C\frac{x}{g(x)}$
    we have, 
    \begin{align*}
      g(x)-g(x-z)
      &=\int_{x-z}^x g'(t) dt 
      \sim z g'(x)(1+o(1))
      %\sim \beta z\frac{g(x)}{x}
      =\beta z\frac{g(x)}{x} + z\frac{g(x)}{x}o(1) \\
      &= \beta z\frac{g(x)}{x} + o(1), 
    \end{align*}  
    as $a\to 0.$ Since $C$ is arbitrary, there exist an increasing function $C(a)\uparrow \infty$ such that 
    for 
    uniformly in  $x\ge (1-\delta)x(a)$ and $z$ such that $z \le C(a)\frac{x}{g(x)}$, 
    \begin{equation}
      \label{def_ca}
      g(x)-g(x-z)=\beta z\frac{g(x)}{x} +o(1). 
    \end{equation}  
    We will split the integral in two parts, 
    \begin{multline*}
      I_2(x,a)
      = 
      2\left(
        \int_{(1-\delta)x(a)}^{x-C(a)\frac{x}{g(x)}}
        + 
        \int_{x-C(a)\frac{x}{g(x)}}^x
      \right)  
    \left(
    \frac{1}{\theta_a}
    +(x-z)
    \right) e^{-\theta_a(x-z)} \overline F(z)dz\\ 
    :=J_1(x,a)+J_2(x,a).
    \end{multline*}  
    %where $\delta(a)\to 0,$ as $a\to 0$. 
    Clearly, uniformly in 
    $x\in [(1-\delta/2)x(a), A_ax(a)]$,  
    \begin{align*}
    J_2(x,a) &= 
    2\int_0^{C(a)\frac{x}{g(x)}} 
    \left(\frac{1}{\theta_a} +z\right) 
    e^{-\theta_a z} \overline F(x-z)\\
    &\sim 
      2x^{-2} 
      \int_0^{C(a)\frac{x}{g(x)}}
      \left(
    \frac{1}{\theta_a}
    +z
    \right) e^{-\theta_a z} e^{-g(x-z)}dz\\ 
    &\sim 
    2x^{-2} e^{-g(x)}
    \int_0^{C(a)\frac{x}{g(x)}}
    \left(
  \frac{1}{\theta_a}
  +z
  \right) e^{-\theta_a z} e^{g(x)-g(x-z)}dz  \\
  &\sim 
  2x^{-2} e^{-g(x)}
  \int_0^{C(a)\frac{x}{g(x)}}
  \left(
\frac{1}{\theta_a}
+z
\right) e^{-\theta_a z} e^{\beta z\frac{g(x)}{x}}dz,
    \end{align*}
 where we used \eqref{def_ca}  in the last step. 
Next note that  uniformly in $x>(1-\delta)x(a)$, 
\begin{align}
  \nonumber
  \theta_a-\beta \frac{g(x)}{x} &\ge
  \theta_a-\beta \frac{g((1-\delta)x(a))}{(1-\delta)x(a)}
  \ge \theta_a - \frac{\beta}{1-\delta} \frac{g(x(a))}{x(a)}\\
  \label{eq.lower.bound.theta}
  &=\theta_a - \frac{\beta}{1-\delta}2a 
  \ge\frac{1}{2} \frac{1-\delta-\beta}{1-\delta}\theta_a,
\end{align}  
as $a\to 0$. 

Changing the variable in  the above integral we obtain 
\begin{align}
  \label{cor4.eq1}
  \frac{x^2e^{g(x)}}{2}
  J_2(x,a)
  &\sim 
 \frac{1}{\theta_a(\theta_a-\beta\frac{g(x)}{x})}
 \int_0^{C(a)\frac{x}{g(x)}(\theta_a-\beta\frac{g(x)}{x})} 
 e^{-u}du\\ 
 \nonumber 
 &\hspace{1cm}+
 \frac{1}{(\theta_a-\beta\frac{g(x)}{x})^2}
 \int_0^{C(a)\frac{x}{g(x)}(\theta_a-\beta\frac{g(x)}{x})} 
 ue^{-u}du
\end{align}  
Now in view of \eqref{eq.lower.bound.theta}, convergence 
$C(a)\to\infty$ and the bound 
$$
\theta_a\frac{x}{g(x)} \ge 
\theta_a\frac{(1-\delta)x(a)}{g((1-\delta)x(a))}
\ge \theta_a(1-\delta) 
\frac{x(a)}{g(x(a))}
=\theta_a(1-\delta)2a \ge (1-2\delta)
$$
both integrals  in \eqref{cor4.eq1}  converges to $1$
uniformly in $x>(1-\delta)x(a)$. Hence, uniformly in $x>(1-\delta)x(a)$, 
\begin{align}
  \nonumber 
  J_2(x,a) &\sim 2x^{-2}e^{-g(x)} 
  \left(
  \frac{1}{\theta_a(\theta_a-\beta\frac{g(x)}{x})}
  +
 \frac{1}{(\theta_a-\beta\frac{g(x)}{x})^2}
 \right)\\ 
 &\sim 
 \label{eq.j2} 
 2\overline F^{I}(x) g'(x)
 \left(
  \frac{1}{\theta_a(\theta_a-\beta\frac{g(x)}{x})}
  +
 \frac{1}{(\theta_a-\beta\frac{g(x)}{x})^2}
 \right)
 .
\end{align}  
To estimate $J_1(x,a)$ we use the 
following estimate:   
for $x\ge (1-\delta)x(a)$, as $a\to 0$, 
for some $\delta_1$ such that $\beta<\delta_1<1-\delta$, 
\begin{align*}
  g(x)-g(z) = \int_z^x g'(t) dt 
  \le \delta_1
  \int_z^x \frac{g(t)}{t}dt 
  \le \delta_1(x-z) \frac{g((1-\delta)x(a))}{(1-\delta)x(a)}\\ 
  \le \frac{\delta_1 }{(1-\delta)} 
  \frac{g(x(a))}{x(a)} (x-z)
  = \frac{\delta_1 2a}{(1-\delta)} 
    (x-z).
\end{align*}
Hence, 
\begin{align*}
  J_1(x,a)
  &\le 4 \int_{(1-\delta)x(a)}^{x-C(a)\frac{x}{g(x)}}(x-z) 
  e^{-\theta_a(x-z)} \overline F(z)dz\\ 
  &\sim 
  4 \int_{(1-\delta)x(a)}^{x-C(a)\frac{x}{g(x)}}(x-z) 
  e^{-\theta_a(x-z)} z^{-2}e^{-g(z)}dz\\ 
  &\le \frac{4}{(1-\delta)^2}
  \frac{1}{x(a)^2} e^{-g(x)}
  \int_{(1-\delta)x(a)}^{x-C(a)\frac{x}{g(x)}}(x-z) 
  e^{-\theta_a(x-z)} e^{g(x)-g(z)}dz
\end{align*}  
Since, $\delta_1/(1-\delta)<1$ and $\theta_a\sim 2a$ 
we obtain for some $\delta_2\in(0,1-\delta_1/(1-\delta)$, 
\begin{align*}
  J_1(x,a) &\le 
  \frac{4}{(1-\delta)^2}
  \frac{1}{x(a)^2} e^{-g(x)}
  \int_{(1-\delta)x(a)}^{x-C(a)\frac{x}{g(x)}}(x-z) 
  e^{-\delta_2 a (x-z)}dz\\ 
  &\le 
  \frac{4}{(1-\delta)^2}
  \frac{1}{x(a)^2} e^{-g(x)}
  \int_{C(a)\frac{x}{g(x)}}^\infty z  
  e^{-\delta_2 a z}dz\\
  &\le 
  \frac{4}{(1-\delta)^2}
  \frac{1}{a^2 x(a)^2} e^{-g(x)}
  \int_{C(a)a\frac{x}{g(x)}}^\infty z  
  e^{-\delta_2 z}dz.%\\0
  %&\le 
  %\frac{4}{(1-\delta)^2\theta_a^2}
  %\frac{1}{x(a)^2} e^{-g(x)} o(1), a\to 0.
\end{align*}  
Note that for $x\ge (1-\delta) x(a)$, the following estimate is valid
$$
\frac{x}{g(x)}\ge 
\frac{(1-\delta) x(a)}{g((1-\delta) x(a))}
\ge \frac{(1-\delta)x(a)}{g(x(a))}
=\frac{(1-\delta)}{2a},
$$
and hence,
$$
C(a)a\frac{x}{g(x)} \ge C(a)(1-\delta)/2\to\infty.
$$
Using the estimate for $t\ge 1$, 
$$
\int_t^\infty z e^{-\delta_2 z}dz \le \frac{2t}{\delta_2^2} e^{-\delta_2 t} 
$$
we obtain with 
 $C=\frac{32}{(1-\delta)^2\delta_2^2}$,  
we 
$$
J_1(x,a)
\le 
C 
\frac{1}{\theta_a^2 x(a)^2} e^{-g(x)}
\frac{C(a)a x}{g(x)} 
\exp\left(
  -\delta_2\frac{C(a)a x}{g(x)} 
\right)
%\int_{C(a)(1-\delta)/2}^\infty z  
%e^{-\delta_2 z}dz.
$$ 
Then,
\begin{align*}
  \frac{J_1(x,a)}{\frac{2}{\theta_a}\overline F^I(x)}
  &\le 
  C
  \frac{
    \frac{1}{\theta_a^2 x(a)^2} e^{-g(x)}
\frac{C(a)a x}{g(x)} 
e^{
  -\delta_2\frac{C(a)a x}{g(x)} 
}
  }
  {\frac{1}{\theta_a}x^{-2}e^{-g(x)}/g'(x)}\\
&= \frac{C}{\theta_a}  
\left(\frac{x}{x(a)}\right)^2
\frac{C(a)a x}{g(x)} 
g'(x)
\exp\left(
  -\delta_2\frac{C(a)a x}{g(x)} 
\right)\\
&\le C 
\left(\frac{x}{x(a)}\right)^2 C(a)
\exp\left(
  -\delta_2\frac{C(a)a x}{g(x)} 
\right)\\
& =C 
\left(\frac{x}{x(a)}\right)^2 C(a)
\exp\left(
  -0.5 \delta_2 C(a) \frac{x}{x(a)}\frac{g(x(a))}{g(x)} 
\right)
.
%  \le 
%  O(1)\widetilde A_a^2 
%  \int_{C(a)(1-\delta)/2}^\infty z  
%    e^{-\delta_2 z}dz
%  \frac{g'(x)}{a}
%  \le o(1)A_a^2 \to 0,
\end{align*}  
Now let $\gamma_0$ be such that \eqref{sc2} holds. Then, 
for $x\ge (1-\delta)x(a)$
\begin{align*}
\frac{x}{x(a)}
\frac{g(x(a))}{g(x)} 
&=
\left(\frac{x}{x(a)}\right)^{1-\gamma_0}
\frac{x^{\gamma_0}}{g(x)}
\frac{g(x(a))}{x(a)^{\gamma_0}}
\ge 
\left(\frac{x}{x(a)}\right)^{1-\gamma_0} 
\frac{((1-\delta)x(a))^{\gamma_0}}{g((1-\delta)x(a))}
\frac{g(x(a))}{x(a)^{\gamma_0}}\\
&\ge (1-\delta)^{\gamma_0} \left(\frac{x}{x(a)}\right)^{1-\gamma_0}.
\end{align*}
Denoting $\delta_3 = 0.5\delta_2 (1-\delta)^{\gamma_0}$ we obtain, 
$$
\frac{J_1(x,a)}{\frac{2}{\theta_a}\overline F^I(x)}
\le 
C 
\left(\frac{x}{x(a)}\right)^2 C(a)
\exp\left(
  -\delta_3 C(a) \left(\frac{x}{x(a)}\right)^{1-\gamma_0} 
\right)
$$
Since $C(a)\to \infty$, $x\ge (1-\delta)x(a)$ 
and $z^2 e^{-z} $ is monotone decreasing to $0$, 
we obtain the required uniform convergence to $0$. 

\end{proof}


\begin{thebibliography}{99}
\bibitem{DFK}
    \textsc{Denisov, D., Foss S.} and \textsc{Korshunov, D.} (2004).
    Tail asymptotics for the supremum of a random walk when the mean is not finite.   
    {\it Queueing Systems } {\bf 46}(1-2), 15–33.
\end{thebibliography}
    \end{document}